\documentclass[11pt]{article}
\usepackage{amssymb}
\usepackage{amsmath}
\usepackage{t1enc}
\usepackage[latin1]{inputenc}
\usepackage[english]{babel}
\usepackage{amsmath,amsthm}
\usepackage{amsfonts}
\usepackage{latexsym}
\usepackage[dvips]{graphicx}
\usepackage{graphicx}
\usepackage[natural]{xcolor}
\usepackage{setspace}

\def\3n{\negthinspace \negthinspace \negthinspace }
\def\2n{\negthinspace \negthinspace }
\def\1n{\negthinspace }

\newtheorem{theorem}{Theorem}[section]
\newtheorem{lemma}{Lemma}[section]
\newtheorem{proposition}{Proposition}[section]

\newtheorem{definition}{Definition}[section]

\newtheorem{remark}{Remark}[section]

\makeatletter
   
   \@addtoreset{equation}{section}
\makeatother
\input{tcilatex}

\begin{document}

\author{Ana Bela Cruzeiro$^{1}\thanks{%
A. B. Cruzeiro acknowledges the financial support of the FCT project
PTDC/MAT/104173/2008.}$, Andr\'{e} de Oliveira Gomes$^{2,3}\thanks{%
A. de Oliveira Gomes acknowledges the financial support partly by the FCT
SFRH/BD/51528/2011 and by the project IRTG 1740/ TRP 2011/50151-0 funded by
the DFG/FAPESP. }$, Liangquan Zhang$^{4,}$$^{5}\thanks{%
L. Zhang (Corresponding author) acknowledges the financial support partly by
the National Nature Science Foundation of China (No. 11201263 and No.
11201264).) and the Nature Science Foundation of Shandong Province (No.
ZR2012AQ004). }$ \\
{\small 1. GFMUL and Dep. de Matemática IST(TUL). }\\
{\small \ Av. Rovisco Pais, 1049-001 Lisboa, Portugal.}\\
{\small \textit{E-mail: abcruz@math.ist.utl.pt.}}\\
{\small 2. GFMUL - Grupo de Física Matemática, Universidade de Lisboa}\\
{\small \ Av. Prof. Gama Pinto 2, 1649-003 Lisboa, Portugal.}\\
{\small 3. Institut für Mathematik Humboldt-Universität zu Berlin, Germany.}%
\\
{\small \textit{E-mail: adeoliveriagomes@sapo.pt.}}\\
{\small 4. School of mathematics, Shandong University, China.}\\
{\small 5. Laboratoire de Mathématiques, }\\
{\small \ \ Université de Bretagne Occidentale, 29285 Brest Cédex, France.}\\
{\small \textit{E-mail: xiaoquan51011@163.com.}}}
\title{Asymptotic Properties of Coupled Forward-Backward Stochastic
Differential Equations}
\maketitle

\begin{abstract}
In this paper, we consider coupled forward-backward stochastic differential
equations (FBSDEs in short) with parameter $\varepsilon >0,$ of the
following type%
\begin{equation*}
\left\{ 
\begin{array}{crl}
X^{\varepsilon ,t,x}\left( s\right) & = & x+\int_{t}^{s}f\left(
r,X^{\varepsilon ,t,x}\left( r\right) ,Y^{\varepsilon ,t,x}\left( r\right)
\right) \mbox{\rm d}r+\sqrt{\varepsilon }\int_{t}^{s}\sigma \left(
r,X^{\varepsilon ,t,x}\left( r\right) ,Y^{\varepsilon ,t,x}\left( r\right)
\right) \mbox{\rm d}W\left( r\right) , \\ 
Y^{\varepsilon ,t,x}\left( s\right) & = & h\left( X^{\varepsilon ,t,x}\left(
T\right) \right) +\int_{s}^{T}g\left( r,X^{\varepsilon ,t,x}\left( r\right)
,Y^{\varepsilon ,t,x}\left( r\right) ,Z^{\varepsilon ,t,x}\left( r\right)
\right) \mbox{\rm d}r \\ 
&  & -\int_{s}^{T}Z^{\varepsilon ,t,x}\left( r\right) \mbox{\rm d}W\left(
r\right) ,\text{ }0\leq t\leq s\leq T.%
\end{array}%
\right.
\end{equation*}%
\noindent We study the asymptotic behavior of its solutions and establish a
large deviation principle for the corresponding processes.
\end{abstract}



\section{Introduction}

Non-linear backward stochastic differential equations (BSDEs in short) were
first introduced in Stochastic Optimal Control Theory with the pioneering
work of Bismut \cite{B} and then developed by Pardoux and Peng \cite{PP1}.
Since then, they have become a powerful tool in many fields, such as
mathematics finance, optimal control, stochastic games, partial differential
equations and homogenization etc. Simultaneously, it is well known that the
Hamiltonian system associated with the maximum principle for stochastic
optimal control problems corresponds to certain fully coupled forward
backward stochastic differential equations (FBSDEs in short) (see \cite{PW}%
). In mathematical finance, fully coupled FBSDEs can be encountered when one
studies the problems of hedging options involved in a large investor in
financial market (see \cite{CM, El}). Besides, fully coupled FBSDEs can
provide probabilistic interpretations for the solutions of a class of
quasilinear parabolic and elliptic PDEs (see \cite{PP2, P1, P2, PW}).

There are three main approaches to solve FBSDEs. The first one, the Method
of Contraction Mapping, which was considered by Antonelli \cite{An} and
later developed by Pardoux-Tang \cite{P2}, works well when the duration $T$
is relatively small. Second, the Four Step Scheme. It removed the
restriction on the time duration for Markovian FBSDEs and was initiated by
Ma-Protter-Yong \cite{MPY}. Third, the Method of Continuation. This method
can treat non-Markovian FBSDEs with arbitrary time duration; it was
initiated by Hu-Peng \cite{HP} and Peng-Wu \cite{PW} and later developed by
Ma et al in \cite{M}. It is worth noting that these three methods do not
cover each other.

In our paper, on one hand, we obtain two kind of asymptotic results. The
first employs the assumptions established by Peng et al in \cite{HP, PW} and
provides the convergence in distribution of the associated processes
(Theorem \ref{t1}). The second, more closely related to the PDE point of
view but still using strictly probabilistic methods, provides almost-sure
convergence with values in some $L^{2}$ type spaces (Theorem \ref{t2}). In
this last result the problem of convergence of classical/viscosity solutions
of the quasilinear parabolic system of PDE's associated to the backward
equation is naturally addressed. Notice that when this PDE takes the form of
the backward Burgers equation, the problem becomes the convergence of
solutions for vanishing viscosity hydrodynamical parameter. We mention the
work \cite{FR} where the authors do the asymptotic studies of FBSDEs with
generalized Burgers type nonlinearities.

The large deviation principle (LDP) characterizes the limiting behavior of
probability measure in applied probability and is largely used in rare
events simulation. Recently, there has been a growing literature on studying
the applications of LDP in finance (see \cite{P3}).

We now consider the following small perturbation of FBSDEs (\ref{1.1})%
\begin{equation}
\left\{ 
\begin{array}{crl}
X^{\varepsilon ,t,x}\left( s\right) & = & x+\int_{t}^{s}f\left(
r,X^{\varepsilon ,t,x}\left( r\right) ,Y^{\varepsilon ,t,x}\left( r\right)
\right) \mbox{\rm d}r \\ 
&  & +\sqrt{\varepsilon }\int_{t}^{s}\sigma \left( r,X^{\varepsilon
,t,x}\left( r\right) ,Y^{\varepsilon ,t,x}\left( r\right) \right) 
\mbox{\rm
d}W\left( r\right) , \\ 
Y^{\varepsilon ,t,x}\left( s\right) & = & h\left( X^{\varepsilon ,t,x}\left(
T\right) \right) +\int_{s}^{T}g\left( r,X^{\varepsilon ,t,x}\left( r\right)
,Y^{\varepsilon ,t,x}\left( r\right) ,Z^{\varepsilon ,t,x}\left( r\right)
\right) \mbox{\rm d}r \\ 
&  & -\int_{s}^{T}Z^{\varepsilon ,t,x}\left( r\right) \mbox{\rm d}W\left(
r\right) ,\text{ }0\leq t\leq s\leq T.%
\end{array}%
\right.  \label{1.1}
\end{equation}%
The solution of this equation is denoted by 
\begin{equation*}
\left( X^{\varepsilon ,t,x}\left( s\right) ,Y^{\varepsilon ,t,x}\left(
s\right) ,Z^{\varepsilon ,t,x}\left( s\right) ,t\leq s\leq T\right) .
\end{equation*}%
We want to establish the large deviation principle of the law of $%
Y^{\varepsilon ,t,x}$ in the space of $C\left( \left[ 0,T\right] ;\mathbb{R}%
^{n}\right) ,$ namely the asymptotic estimates of probabilities $P\left(
Y^{\varepsilon ,t,x}\in \Gamma \right) ,$ where $\Gamma \in \mathcal{B}%
\left( C\left( \left[ 0,T\right] ;\mathbb{R}^{n}\right) \right) .$

Ma et al in \cite{MaZ} first considered the sample path large deviation
principle for the adapted solutions to the FBSDEs in the case where the
drift term contains $Z$ and under appropriate conditions in terms of a
certain type of convergence of solutions for the associated quasilinear
PDEs. While, with probability methods, by the contraction principle, the
same small random perturbation for BSDEs and the Freidlin-Wentzell's large
deviation estimates in $C\left( \left[ 0,T\right] ;\mathbb{R}^{n}\right) $
are also obtained by \cite{D, Es, FR, R, KM} (see references therein). In 
\cite{CX} a large deviation principle of Freidlin and Wentzell type under
nonlinear probability for diffusion processes with a small diffusion
coefficient was obtained.

Our aim in this paper extends the previous work \cite{D, Es, FR, R, KM} to
the coupled case. Namely, when the drift and diffusion terms contain $Y$, by
probability methods and under some suitable assumptions. Note that in \cite%
{CX}, the authors considered the fully coupled FBSDEs via the corresponding
PDE techniques.

Our work is organized as follows. In Section \ref{s1}, we give the framework
of our paper. Some estimates and regularity results are established for the
solutions of FBSDEs (\ref{1.1}) in Section \ref{s1}. Then in Section \ref{s2}
we show our main results Theorem \ref{t1} and Theorem \ref{t2}. Section \ref%
{s3} is devoted to establishing the large deviation results for (\ref{1.1}).
Some proofs of technique lemmas are given in the Appendix.

\section{Preliminaries}

\label{s1}

Let us begin by introducing the setting for the stochastic differential
systems we want to investigate. Consider as Brownian motion $W$ the $d$%
-dimensional coordinate process on the classical Wiener space $\left( \Omega
,\mathcal{F},P\right) $, i.e., $\Omega $ is the set of continuous functions
from $\left[ 0,T\right] $ to $\mathbb{R}^{d}$ starting from $0$, $\Omega
=C\left( \left[ 0,T\right] ;\mathbb{R}^{d},0\right) $, $\mathcal{F}$ the
completed Borel $\sigma $-algebra over $\Omega $, $P$ the Wiener measure and 
$W$ the canonical process: $W_{s}\left( \omega \right) =\omega _{s}$, $s\in %
\left[ 0,T\right] ,$ $\omega \in \Omega .$ By $\left\{ \mathcal{F}_{s},0\leq
s<T\right\} $ we denote the natural filtration generated by $\left\{
W_{s}\right\} _{0\leq s<T}$ and augmented by all $P$-null sets, i.e., 
\begin{equation*}
\mathcal{F}_{s}=\sigma \left\{ W\left( r\right) ,r\leq s\right\} \vee 
\mathcal{N}_{p},\text{ }s\in \left[ 0,T\right] ,
\end{equation*}%
where $\mathcal{N}_{p}$ is the set of all $P$-null subsets. For each $t>0,$
we denote by $\left\{ \mathcal{F}_{s}^{t},\text{ }t\leq s\leq T\right\} $
the natural filtration of the Brownian motion $\left\{ W\left( s\right)
-W\left( t\right) ,\text{ }t\leq s\leq T\right\} ,$ augmented by the $P$%
-null set of $\mathcal{F}$.

We denote by $\mathcal{M}^{2}\left( t,T;\mathbb{R}^{p}\right) $ the set of
all $\mathbb{R}^{p}$-valued $\mathcal{F}_{t}$-adapted process $v\left( \cdot
\right) $ such that 
\begin{equation*}
\mathbb{E}\left[ \int_{t}^{T}\left\vert v\left( s\right) \right\vert ^{2}%
\text{d}s\right] <+\infty ,
\end{equation*}%
and $\mathcal{S}^{2}\left( t,T;\mathbb{R}^{n}\right) ,$ the set of all $%
\mathbb{R}^{n}$-valued $\mathcal{F}_{t}$-adapted process $v\left( \cdot
\right) $ such that%
\begin{equation*}
\mathbb{E}\left[ \sup\limits_{t\leq s\leq T}\left\vert v\left( s\right)
\right\vert ^{2}\right] <+\infty .
\end{equation*}%
Clearly, $\mathcal{N}_{t}=\mathcal{S}^{2}\left( t,T;\mathbb{R}^{n}\right)
\times \mathcal{S}^{2}\left( t,T;\mathbb{R}^{n}\right) \times \mathcal{M}%
^{2}\left( t,T;\mathbb{R}^{p}\right) $ forms a Banach space.

The space $\mathcal{N}_{t}$ is a Banach space for the natural norm
associated with the product topology structure.

We use the usual inner product and Euclidean norm in $\mathbb{R}^{n}$, $%
\mathbb{R}^{m},$ and $\mathbb{R}^{m\times d}.$

\begin{definition}
\label{d1}A triple of processes 
\begin{equation*}
\left( X^{\varepsilon ,t,x},Y^{\varepsilon ,t,x},Z^{\varepsilon ,t,x}\right)
:\left[ 0,T\right] \times \Omega \rightarrow \mathbb{R}^{n}\times \mathbb{R}%
^{m}\times \mathbb{R}^{m\times d}
\end{equation*}%
is called an adapted solution of the Eqs. (\ref{1.1}), if $\left(
X^{\varepsilon ,t,x},Y^{\varepsilon ,t,x},Z^{\varepsilon ,t,x}\right) \in 
\mathcal{N}_{t},$ and it satisfies (\ref{1.1}) $P$-a.s..
\end{definition}

It is clear that the above Eq. (\ref{1.1}) is a stochastic two point
boundary value problem. Especially, it contains a deterministic two point
boundary value problem as a special case when $\varepsilon \rightarrow 0$.

We are given an $m\times n$ full-rank matrix $G$. We use the notations%
\begin{equation*}
u=\left( 
\begin{array}{c}
x \\ 
y \\ 
z%
\end{array}%
\right) ,\qquad A^{\varepsilon }\left( t,u\right) =\left( 
\begin{array}{c}
-G^{T}g \\ 
Gf \\ 
\sqrt{\varepsilon }G\sigma%
\end{array}%
\right) \left( t,u\right) .
\end{equation*}%
We now give the first assumptions of our paper:

\begin{enumerate}
\item[\textbf{(A1)}] Let $m\geq n.$ The processes $g\left( \cdot
,x,y,z\right) ,$ $f\left( \cdot ,x,y\right) $ and $\sigma \left( \cdot
,x,y\right) $ are $\mathcal{F}_{t}$-adapted, and the random variable $%
h\left( x\right) $ is $\mathcal{F}_{T}$-measurable, for all $\left(
x,y,z\right) \in \mathbb{R}^{n}\times \mathbb{R}^{m}\times \mathbb{R}%
^{m\times d}$. $f,$ $\sigma $ and $g$ are continuous in $t,$ $P$-a.s.$.$
Moreover, the following holds:%
\begin{equation*}
\mathbb{E}\left[ \left\vert h\left( 0\right) \right\vert ^{2}\right]
<+\infty .
\end{equation*}

There exists a constant $C_{1}>0,$ such that, $u^{i}=\left(
x^{i},y^{i},z^{i}\right) \in \mathbb{R}^{n}\times \mathbb{R}^{m}\times 
\mathbb{R}^{m\times d},$ $i=1,$ $2,$%
\begin{equation*}
\begin{array}{rcl}
\left\vert g\left( t,x^{1},y^{1},z^{1}\right) -g\left(
t,x^{2},y^{2},z^{2}\right) \right\vert & \leq & C_{1}\left( \left\vert
x^{1}-x^{2}\right\vert +\left\vert y^{1}-y^{2}\right\vert +\left\vert
z^{1}-z^{2}\right\vert \right) , \\ 
\left\vert f\left( t,x^{1},y^{1}\right) -f\left( t,x^{2},y^{2}\right)
\right\vert & \leq & C_{1}\left( \left\vert x^{1}-x^{2}\right\vert
+\left\vert y^{1}-y^{2}\right\vert \right) , \\ 
\left\vert \sigma \left( t,x^{1},y^{1}\right) -\sigma \left(
t,x^{2},y^{2}\right) \right\vert & \leq & C_{1}\left( \left\vert
x^{1}-x^{2}\right\vert +\left\vert y^{1}-y^{2}\right\vert \right) , \\ 
\left\vert h\left( x^{1}\right) -h\left( x^{2}\right) \right\vert & \leq & 
C_{1}\left\vert x^{1}-x^{2}\right\vert , \\ 
&  & P\text{-a.s., a.e. }t\in \mathbb{R}^{+}.%
\end{array}%
\end{equation*}

\item[(\textbf{A2)}] There exists a constant $C_{2}>0,$ such that for any
fixed $\varepsilon \in \left( 0,1\right] ,$%
\begin{eqnarray*}
&&\left\langle A^{\varepsilon }\left( t,u^{1}\right) -A^{\varepsilon }\left(
t,u^{2}\right) ,u^{1}-u^{2}\right\rangle \\
&\leq &-\left( C_{2}+\sqrt{\varepsilon }\right) \left( \left\vert G\left(
x^{1}-x^{2}\right) \right\vert ^{2}+\left\vert G^{T}\left(
y^{1}-y^{2}\right) \right\vert ^{2}\right) ,\text{ }P\text{-a.s., a.e. }t\in 
\mathbb{R}^{+},
\end{eqnarray*}

and%
\begin{equation*}
\left\langle h\left( x^{1}\right) -h\left( x^{2}\right)
,x^{1}-x^{2}\right\rangle \geq C_{2}\left\vert G\left( x^{1}-x^{2}\right)
\right\vert ^{2},\text{ }P\text{-a.s., }\forall \left( x^{1},x^{2}\right)
\in \mathbb{R}^{n}\times \mathbb{R}^{n}.
\end{equation*}
\end{enumerate}

\noindent We have the following:

\begin{proposition}
\label{pr1}Assume that (A1) and (A2) hold, then there exists a unique
adapted solution 
\begin{equation*}
\left( X^{\varepsilon ,t,x},Y^{\varepsilon ,t,x},Z^{\varepsilon ,t,x}\right)
\end{equation*}
for Eqs. (\ref{1.1}).
\end{proposition}

\noindent The proof can be seen in Theorem 2.6, Remark 2.8 in \cite{PW}.

\begin{remark}
\label{re1}\textsl{If }$m=n,$ \textsl{then $G=I_{n}.$} \textsl{\ If }$m<n,$%
\textsl{\ according to Theorem 2.6 in \cite{PW}, the first condition of (A2)
should be }%
\begin{eqnarray*}
&&\left\langle A^{\varepsilon }\left( t,u^{1}\right) -A^{\varepsilon }\left(
t,u^{2}\right) ,u^{1}-u^{2}\right\rangle \\
&\leq &-\left( C_{2}+\sqrt{\varepsilon }\right) \left( \left\vert G\left(
z^{1}-z^{2}\right) \right\vert ^{2}+\left\vert G^{T}\left(
y^{1}-y^{2}\right) \right\vert ^{2}\right) ,\text{ }P\text{-a.s., a.e. }t\in 
\mathbb{R}^{+}.
\end{eqnarray*}%
\textsl{However, if setting }$x^{1}=x^{2},$\textsl{\ }$y^{1}=y^{2},$\textsl{%
\ there exists a contradiction, because }$\sigma $\textsl{\ does not contain 
}$z$\textsl{.}
\end{remark}

\noindent From now on, let us suppose that 
\begin{equation*}
\text{the coefficients }f\text{, }\sigma \text{, }g\text{ and }h\text{ are
deterministic, i.e., independent of }\omega \in \Omega .
\end{equation*}%
Now consider 
\begin{equation}
u^{\varepsilon }\left( t,x\right) =Y^{\varepsilon ,t,x}\left( t\right)
,\qquad \left( t,x\right) \in \left[ 0,T\right] \times \mathbb{R},\text{ }%
\varepsilon \in \left( 0,1\right] ,  \label{2.1}
\end{equation}%
which is a deterministic vector since it is $\mathcal{F}_{t}^{t}$ measurable
(Blumenthal's \textit{0-1 Law} -see Remark 1.2 of \cite{De}). In \cite{P1}
it is shown that $u$ is a viscosity solution of the associated quasilinear
parabolic partial differential equation 
\begin{equation}
\begin{cases}
\dfrac{\partial (u^{\varepsilon })^{l}}{\partial t}(t,x)+g^{l}(t,x,u^{%
\varepsilon }(t,x),\sqrt{\varepsilon }\nabla _{x}u^{\varepsilon }(t,x)\sigma
^{\varepsilon }(t,x,u^{\varepsilon }(t,x))) \\ 
\qquad +\sum_{i=1}^{n}f_{i}^{l}(t,x,u^{\varepsilon }(t,x))\dfrac{\partial
(u^{\varepsilon })^{l}}{\partial x_{i}}(t,x)+\dfrac{\varepsilon }{2}%
\sum_{i=1}^{n}a_{ij}(t,x,u^{\varepsilon }(t,x))\dfrac{\partial
^{2}(u^{\varepsilon })^{l}}{\partial x_{i}\partial x_{j}}(t,x)=0, \\ 
u^{\varepsilon }(T,x)=h(x),\,x\in \mathbb{R}^{n};\,\,t\in \lbrack
0,T];\,\,l=1,...,n,%
\end{cases}
\label{2.2}
\end{equation}%
where $a_{i,j}=(\sigma \sigma ^{T})_{i,j}$.

We define below the notion of viscosity solution for the parabolic system of
Partial Differential Equations (PDEs for short) (\ref{2.2}). For each $%
\varepsilon >0$, consider the following differential operator, 
\begin{align*}
& (L_{l}^{\varepsilon }\varphi )(t,x,y,z)=\dfrac{\varepsilon }{2}%
\displaystyle\sum_{i,j=1}^{n}a_{ij}(t,x,y)\dfrac{\partial ^{2}\varphi ^{l}}{%
\partial x_{i}\partial x_{j}}(t,x)+<f(t,x,y),\nabla \varphi ^{l}(t,x)> \\
& l=1,...,n,\text{ }\forall \varphi \in C^{1,2}([0,T]\times \mathbb{R}^{n},%
\mathbb{R}^{n}),\,t\in \lbrack 0,T],\,(x,y,z)\in \mathbb{R}^{n}\times 
\mathbb{R}^{n}\times \mathbb{R}^{n\times d}.
\end{align*}

\noindent The space $C^{1,2}([0,T]\times \mathbb{R}^{n},\mathbb{R}^{n})$ is
the space of the functions $\phi :[0,T]\times \mathbb{R}^{n}\rightarrow 
\mathbb{R}^{n}$, which are $C^{1}$ with respect to the first variable and $%
C^{2}$ with respect to the second variable.

The system (\ref{2.2}) reads%
\begin{equation}
\left\{ 
\begin{array}{l}
\dfrac{\partial (u^{\varepsilon })^{l}}{\partial t}+(L_{l}^{\varepsilon
}u^{\varepsilon })(t,x,u^{\varepsilon }(t,x),\nabla _{x}u^{\varepsilon
}(t,x)\sigma (t,x,u^{\varepsilon }(t,x)))+ \\ 
\qquad g^{l}(t,x,u^{\varepsilon }(t,x),\sqrt{\varepsilon }\nabla
_{x}u^{\varepsilon }(t,x)\sigma (t,x,u^{\varepsilon }(t,x)))=0, \\ 
u^{\varepsilon }(T,x)=h(x),\,x\in \mathbb{R}^{d},\,\,\,t\in \lbrack 0,T],%
\text{ }l=1,...,n.%
\end{array}%
\right.  \label{2.3}
\end{equation}

\begin{definition}[\textbf{Viscosity Solutions}]
\label{de2}Let $u^{\varepsilon }\in C\big (\lbrack 0,T]\times \mathbb{R}^{n},%
\mathbb{R}^{n}\big )$. The function $u^{\varepsilon }$ is said to be a
viscosity sub-solution (resp. super-solution) of the system (\ref{2.2}) if%
\begin{equation*}
(u^{\varepsilon })^{l}(T,x)\leq h^{l}(x);\,\,\forall l=1,...,n;\,\,x\in 
\mathbb{R}^{n}
\end{equation*}%
(resp. $(u^{\varepsilon })^{l}(T,x)\geq h^{l}(x);\,\,\forall
l=1,...,n;\,\,x\in \mathbb{R}^{n}$ )

and for each $l=1,...,n,(t,x)\in \lbrack 0,T]\times \mathbb{R}^{n},$ $%
\varphi \in C^{1,2}([0,T]\times \mathbb{R}^{n},\mathbb{R}^{n})$ such that $%
(t,x)$ is a local minimum (resp. maximum) point of $\varphi -u^{\varepsilon
} $, we have%
\begin{equation*}
\begin{array}{l}
\dfrac{\partial \varphi }{\partial t}(t,x)+(L_{l}^{\varepsilon }\varphi
)(t,x,u^{\varepsilon }(t,x),\sqrt{\varepsilon }\nabla _{x}\varphi
(t,x)\sigma (t,x,u^{\varepsilon }(t,x)))+ \\ 
g^{l}(t,x,u^{\varepsilon }(t,x),\sqrt{\varepsilon }\nabla _{x}\varphi
(t,x)\sigma (t,x,u^{\varepsilon }(t,x)))\geq 0;l=1,...,n%
\end{array}%
\end{equation*}%
(resp. $\leq 0).$

The function $u^{\varepsilon }$ is said to be a viscosity solution of the
system (\ref{2.2}) if $u^{\varepsilon }$ is both a viscosity sub-solution
and a viscosity super-solution of this system.
\end{definition}

Under more restrictive assumptions (that we present in Section \ref{s2}), we
shall have 
\begin{equation}
u^{\varepsilon }(t,x)=Y^{t,\varepsilon ,x}\left( t\right)  \label{2.4}
\end{equation}%
and $u^{\varepsilon }$ will be actually a classical solution of (\ref{2.2}).
This can be proved using the \textit{Four Step Scheme Methodology} of
Ma-Protter-Yong \cite{MPY} with the help of Ladyzhenskaja's work in
quasilinear parabolic PDEs \cite{LS}.

For simplicity, we only consider the case where both $X^{\varepsilon ,t,x}$
and $Y^{\varepsilon ,t,x}$ are $n$-dimensional, that is $m=n$, then $G=I_{n}$%
. The result is analogous for the case $m>n.$

Next we introduce another set of assumptions which is slightly different
from (A1)-(A2).

\begin{enumerate}
\item[(\textbf{A3)}] We say $f,$ $g,$ $\sigma $ and $h$ satisfy (A3) if
there exist two constants $C_{1},$ $\Lambda >0$ such that: (A1) hold as well
as:

\noindent \textit{(A.3.1)} $\forall t\in \lbrack 0,T]$, $\forall
(x_{1},y_{1},z_{1}),$ $(x_{2},y_{2},z_{2})\in \mathbb{R}^{n}\times \mathbb{R}%
^{n}\times \mathbb{R}^{n\times d}$:%
\begin{eqnarray*}
\left\langle x_{1}-x_{2},f(t,x_{1},y_{1})-f(t,x_{2},y_{1})\right\rangle
&\leq &C_{1}\left\vert x_{1}-x_{2}\right\vert ^{2}, \\
\left\langle y_{1}-y_{2},f(t,x_{1},y_{1})-f(t,x,_{1}y_{2})\right\rangle
&\leq &C_{1}\left\vert y_{1}-y_{2}\right\vert ^{2}.
\end{eqnarray*}%
\textit{(A.3.2)} $\forall \,t\in \lbrack 0,T]$, $\forall (x,y,z)\in \mathbb{R%
}^{n}\times \mathbb{R}^{n}\times \mathbb{R}^{n\times d}$:\vspace{1pt}%
\begin{eqnarray*}
\left\vert f(t,x,y)\right\vert &\leq &\Lambda (1+\left\vert x\right\vert
+\left\vert y\right\vert ), \\
\left\vert g(t,x,y,z)\right\vert &\leq &\Lambda (1+\left\vert x\right\vert
+\left\vert y\right\vert +\left\vert z\right\vert ), \\
\left\vert \sigma (t,x,y)\right\vert &\leq &\Lambda (1+\left\vert
x\right\vert +\left\vert y\right\vert ), \\
\left\vert h(x)\right\vert &\leq &\Lambda (1+\left\vert x\right\vert ).
\end{eqnarray*}%
\textit{(A.3.3)} $\forall \,t\in \lbrack 0,T]$, $\forall (x,y,z)\in \mathbb{R%
}^{n}\times \mathbb{R}^{n}\times \mathbb{R}^{n\times d}$:%
\begin{eqnarray*}
u &\mapsto &f(t,u,y), \\
v &\mapsto &g(t,x,v,z)\,\text{are\thinspace continuous\thinspace mappings}.
\end{eqnarray*}
\end{enumerate}

\noindent Under this set of hypothesis, Theorem 1.1 of \cite{De} ensures
that there exists a constant $C=C(L)>0$ , depending only on $C_{1}$, such
that for every $T\leq C,$ (\ref{1.1}) admits a unique solution in $\mathcal{N%
}_{t}$.

Moreover, using Theorem 2.6 of \cite{De}, we have 
\begin{equation}
\mathbb{P}\left( \forall \,s\in \lbrack t,T]:u^{\varepsilon }(s,X_{s}^{t,\xi
})=Y_{s}^{t,\xi }\right) =1.  \label{2.5}
\end{equation}%
\begin{equation}
\mathbb{P}\otimes \mu \big (\big \{(\omega ,s)\in \Omega \times \lbrack
t,T]:\left\vert Z_{s}^{t,\varepsilon ,x}(\omega )\right\vert \geq \Gamma _{1}%
\big \}\big )=0.  \label{2.6}
\end{equation}%
where $\mu $ stands for the Lebesgue measure in the real line and $\Gamma
_{1}$ is a constant which only depends on $C_{1},$ $\Lambda ,$ $n,$ $d,$ $T$.

By Remark 2.7 of \cite{De} we have that, for each $\varepsilon >0$, $%
u^{\varepsilon }$ only depends on the coefficients of the system (\ref{2.1}%
), $f,$ $g,$ $\sqrt{\varepsilon }\sigma ,$ $h$. The fact that the dependence
of $\Gamma _{1}=\Gamma _{1}(C_{1},\Lambda ,n,T)$ determines that the
properties (\ref{2.5}), (\ref{2.6}) above hold uniformly in $\varepsilon $;
in particular, there exist continuous versions of $(Y_{s}^{t,\varepsilon
,x},Z_{s}^{t,\varepsilon ,x})_{t\leq s\leq T}$ which are uniformly bounded,
also in $\varepsilon $.

We now assume more regularity on the coefficients of (\ref{1.1}) in the next
set of assumptions.

\begin{enumerate}
\item[(\textbf{A4)}] For $T\leq C$, we say that $f,$ $g,$ $h,$ $\sigma $
satisfy (A4) if there exists three constants $\lambda $, $\Lambda ,$ $\gamma
>0$

\noindent \textit{(A.4.1)} $\forall t\in \lbrack 0,T]$, $\forall (x,y,z)\in 
\mathbb{R}^{n}\times \mathbb{R}^{n}\times \mathbb{R}^{n\times d}$:%
\begin{eqnarray*}
\left\vert f(t,x,y)\right\vert &\leq &\Lambda (1+\left\vert y\right\vert ),
\\
\left\vert g(t,x,y,z)\right\vert &\leq &\Lambda \left( 1+\left\vert
y\right\vert +\left\vert z\right\vert \right) , \\
\left\vert \sigma (t,x,y)\right\vert &\leq &\Lambda , \\
\left\vert h(x)\right\vert &\leq &\Lambda.
\end{eqnarray*}

\textit{(A.4.2) }$\forall (t,x,y)\in \lbrack 0,T]\times \mathbb{R}^{n}\times 
\mathbb{R}^{n}\times \mathbb{R}^{n\times d}:\left\langle \xi ,a(t,x,y)\xi
\right\rangle \geq \lambda \left\vert \xi \right\vert ^{2},\,\forall \xi \in 
\mathbb{R}^{n},$ where $a(t,x,y)=\sigma \sigma ^{T}(t,x,y).$

\textit{(A.4.3)} The function $\sigma $ is continuous.

\textit{(A.4.4)} The function $\sigma $ is differentiable with respect to $x$
and $y$ and its derivatives with respect to $x$ and $y$ are $\gamma $- Hö%
lder in $x$ and $y$, uniformly in $t$.
\end{enumerate}

Under the set of assumptions (A.3) and (A.4), using Proposition 2.4 and
Proposition B.6 of \cite{De} and Theorem 2.9 of \cite{Dee}, one can prove
that there exist two constants $\kappa ,$ $\kappa _{1}$, only depending on $%
L,$ $\Lambda ,$ $n,$ $T$ (independent of $\varepsilon $) such that: 
\begin{equation}
\left\vert u^{\varepsilon }(t,x)\right\vert \leq \kappa ,  \label{2.7}
\end{equation}%
\begin{equation}
u^{\varepsilon }\in C_{b}^{1,2}([0,T]\times \mathbb{R}^{n}),  \label{2.8}
\end{equation}%
\begin{equation}
\sup_{(t,x)\in \lbrack 0,T]\times \mathbb{R}^{d}}\left\vert \nabla
_{x}u^{\varepsilon }(t,x)\right\vert \leq \kappa _{1},  \label{2.9}
\end{equation}%
and

\begin{equation}
Z_{s}^{t,\varepsilon ,x}=\sqrt{\varepsilon }\nabla _{x}u^{\varepsilon
}(t,X_{s}^{t,\varepsilon ,x})\sigma (s,X_{s}^{t,\varepsilon
},Y_{s}^{t,\varepsilon ,x}).  \label{2.10}
\end{equation}%
where $u^{\varepsilon }$ solves uniquely (\ref{2.2}) in $C_{b}^{1,2}([0,T]%
\times \mathbb{R}^{n})$, the space of $C^{1}$ functions with respect to the
first variable and $C^{2}$ with respect to the second variable, with bounded
derivatives. All these facts can be proven probabilistically. The last claim
and the properties (\ref{2.7}), (\ref{2.8}) and (\ref{2.10}) are proved in 
\cite{De}. Delarue delivers in the appendices of \cite{De} probabilistic
methods to obtain these regularity results, under assumptions (A.3) and
(A.4) and over a small time enough duration, using \textit{Malliavin}
Calculus techniques. The estimate of the gradient (\ref{2.9}) is established
by a probabilistic scheme in \cite{Dee}, using a variant of the \textit{%
Malliavin-Bismut} integration by parts formula proposed by Thalmaier \cite{T}
and applied in \cite{TW} to establish a gradient estimate of interior type
for the solutions of a linear elliptic equation on a manifold.

To start with, let us fix $x,$ $y\in \mathbb{R}^{n}$, $\varepsilon \in
\left( 0,1\right] .$ We establish second order moment estimates for the
solution of FBSDEs (\ref{1.1}), which will be essential in Section \ref{s2}.
For convenience, we use the following notations in this paper:%
\begin{equation*}
\left\{ 
\begin{array}{crl}
f^{\varepsilon ,t,x}\left( r\right) & = & f\left( r,X^{\varepsilon
,t,x}\left( r\right) ,Y^{\varepsilon ,t,x}\left( r\right) \right) , \\ 
\sigma ^{\varepsilon ,t,x}\left( r\right) & = & \sigma \left(
r,X^{\varepsilon ,t,x}\left( r\right) ,Y^{\varepsilon ,t,x}\left( r\right)
\right) , \\ 
g^{\varepsilon ,t,x}\left( r\right) & = & g\left( r,X^{\varepsilon
,t,x}\left( r\right) ,Y^{\varepsilon ,t,x}\left( r\right) ,Z^{\varepsilon
,t,x}\left( r\right) \right) , \\ 
h^{\varepsilon ,t,x}\left( T\right) & = & h\left( X^{\varepsilon ,t,x}\left(
T\right) \right) , \\ 
& \vdots &  \\ 
& \text{etc.} & 
\end{array}%
\right.
\end{equation*}%
The proof of the following results can be found in the Appendix.

\begin{lemma}
\label{l1}Assume that (A1) and (A2) hold. Then we have%
\begin{equation}
\left\{ 
\begin{array}{ccc}
\mathbb{E}\left[ \sup\limits_{t\leq s\leq T}\left\vert X^{\varepsilon
,t,x}\left( s\right) -X^{\varepsilon ,t,y}\left( s\right) \right\vert ^{2}%
\right] & \leq & \mathcal{C}_{1}\left\vert x-y\right\vert ^{2}, \\ 
\mathbb{E}\left[ \sup\limits_{t\leq s\leq T}\left\vert Y^{\varepsilon
,t,x}\left( s\right) -Y^{\varepsilon ,t,y}\left( s\right) \right\vert ^{2}%
\right] & \leq & \mathcal{C}_{1}\left\vert x-y\right\vert ^{2}, \\ 
\mathbb{E}\left[ \int_{t}^{T}\left\vert Z^{\varepsilon ,t,x}\left( s\right)
-Z^{\varepsilon ,t,y}\left( s\right) \right\vert ^{2}\mbox{\rm d}s\right] & 
\leq & \mathcal{C}_{1}\left\vert x-y\right\vert ^{2},%
\end{array}%
\right.  \label{2.11}
\end{equation}%
where $\mathcal{C}_{1}$ is a positive constant independent of $\varepsilon $
and $t.$
\end{lemma}

\begin{lemma}
\label{l2}Assume that (A1) and (A2) hold. Then we have%
\begin{equation}
\left\{ 
\begin{array}{ccc}
\mathbb{E}\left[ \sup\limits_{t\leq s\leq T}\left\vert X^{\varepsilon
,t,x}\left( s\right) \right\vert ^{2}\right] & \leq & \mathcal{C}_{2}\left(
1+\left\vert x\right\vert ^{2}\right) , \\ 
\mathbb{E}\left[ \sup\limits_{t\leq s\leq T}\left\vert Y^{\varepsilon
,t,x}\left( s\right) \right\vert ^{2}\right] & \leq & \mathcal{C}_{2}\left(
1+\left\vert x\right\vert ^{2}\right) , \\ 
\mathbb{E}\left[ \int_{t}^{T}\left\vert Z^{\varepsilon ,t,x}\left( s\right)
\right\vert ^{2}\mbox{\rm d}s\right] & \leq & \mathcal{C}_{2}\left(
1+\left\vert x\right\vert ^{2}\right) ,%
\end{array}%
\right.  \label{2.12}
\end{equation}%
where $\mathcal{C}_{2}$ is a positive constant independent of $\varepsilon $
and $t.$
\end{lemma}

The following Lemma shows the continuity on $t.$

\begin{lemma}
\label{l3}Assume that (A1) and (A2) hold, then we have.%
\begin{equation}
\left\{ 
\begin{array}{lll}
\mathbb{E}\left[ \sup\limits_{t_{1}\vee t_{2}\leq s\leq T}\left\vert
X^{\varepsilon ,t_{1},x}\left( s\right) -X^{\varepsilon ,t_{2},x}\left(
s\right) \right\vert ^{2}\right] & \leq & \mathcal{C}_{3}\left\vert
t_{1}-t_{2}\right\vert \left( 1+\left\vert x\right\vert ^{2}\right) , \\ 
\mathbb{E}\left[ \sup\limits_{t_{1}\vee t_{2}\leq s\leq T}\left\vert
Y^{\varepsilon ,t_{1},x}\left( s\right) -Y^{\varepsilon ,t_{2},x}\left(
s\right) \right\vert ^{2}\right] & \leq & \mathcal{C}_{3}\left\vert
t_{1}-t_{2}\right\vert \left( 1+\left\vert x\right\vert ^{2}\right) , \\ 
\mathbb{E}\left[ \int_{t_{1}}^{T}\left\vert Z^{\varepsilon ,t_{1},x}\left(
s\right) -Z^{\varepsilon ,t_{2},x}\left( s\right) \right\vert ^{2}%
\mbox{\rm
d}s\right] & \leq & \mathcal{C}_{3}\left\vert t_{1}-t_{2}\right\vert \left(
1+\left\vert x\right\vert ^{2}\right) ,%
\end{array}%
\right.  \label{2.13}
\end{equation}%
where $\mathcal{C}_{4}$ is a positive constant independent of $\varepsilon .$
\end{lemma}

\begin{lemma}
\label{l4}Assume that (A1) and (A2) hold. Pick $0<\varepsilon
_{2}<\varepsilon _{1}<1.$ Then we have%
\begin{equation}
\left\{ 
\begin{array}{ccc}
\mathbb{E}\left[ \sup\limits_{t\leq s\leq T}\left\vert X^{\varepsilon
_{1},t,x}\left( s\right) -X^{\varepsilon _{2},t,x}\left( s\right)
\right\vert ^{2}\right] & \leq & \mathcal{C}_{4}\left( \sqrt{\varepsilon _{1}%
}-\sqrt{\varepsilon _{2}}\right) , \\ 
\mathbb{E}\left[ \sup\limits_{t\leq s\leq T}\left\vert Y^{\varepsilon
_{1},t,x}\left( s\right) -Y^{\varepsilon _{2},t,x}\left( s\right)
\right\vert ^{2}\right] & \leq & \mathcal{C}_{4}\left( \sqrt{\varepsilon _{1}%
}-\sqrt{\varepsilon _{2}}\right) , \\ 
\mathbb{E}\left[ \int_{t}^{T}\left\vert Z^{\varepsilon _{1},t,x}\left(
s\right) -Z^{\varepsilon _{2},t,x}\left( s\right) \right\vert ^{2}%
\mbox{\rm
d}s\right] & \leq & \mathcal{C}_{4}\left( \sqrt{\varepsilon _{1}}-\sqrt{%
\varepsilon _{2}}\right) ,%
\end{array}%
\right.  \label{2.14}
\end{equation}%
where $\mathcal{C}_{3}$ is a positive constant independent of $\varepsilon $
and $t.$
\end{lemma}

Now consider the following deterministic equations 
\begin{equation}
\left\{ 
\begin{array}{rcl}
\mathcal{X}^{t,x}\left( s\right) & = & x+\int_{t}^{s}f\left( r,\mathcal{X}%
^{t,x}\left( r\right) ,\mathcal{Y}^{t,x}\left( r\right) ,0\right) 
\mbox{\rm
d}r, \\ 
\mathcal{Y}^{t,x}\left( s\right) & = & h\left( \mathcal{X}^{t,x}\left(
T\right) \right) +\int_{s}^{T}g\left( r,\mathcal{X}^{t,x}\left( r\right) ,%
\mathcal{Y}^{t,x}\left( r\right) ,0\right) \mbox{\rm d}r.%
\end{array}%
\right.  \label{2.15}
\end{equation}%
By the regularity of $f,$ $g,$ $h$ we have the following

\begin{lemma}
\label{l5}Assume that (A1) and (A2) hold, then there exists a unique
solution $\left( \mathcal{X}^{t,x},\mathcal{Y}^{t,x}\right) $ for Eqs.
(2.15).
\end{lemma}

Clearly, the estimates in Lemma \ref{l4} shows that the triple of 
\begin{equation*}
\left( X^{\varepsilon ,t,x}\left( s\right) ,Y^{\varepsilon ,t,x}\left(
s\right) ,Z^{\varepsilon ,t,x}\left( s\right) \right) _{s\in \left[ t,T%
\right] }
\end{equation*}%
are Cauchy sequences and therefore converge in 
\begin{equation*}
\left( \mathcal{S}^{2}\left( t,T;\mathbb{R}\right) \times \mathcal{S}%
^{2}\left( t,T;\mathbb{R}\right) \times \mathcal{M}^{2}\left( t,T;\mathbb{R}%
^{d}\right) \right) .
\end{equation*}%
Denote the limit by $\left( \mathcal{X}^{0,t,x}\left( s\right) ,\mathcal{Y}%
^{0,t,x}\left( s\right) ,0\right) _{s\in \left[ t,T\right] },$ as $%
\varepsilon \rightarrow 0.$ By uniqueness and existence of Eqs. (\ref{2.15}%
), we know that the limit $\left( \mathcal{X}^{0,t,x}\left( \cdot \right) ,%
\mathcal{Y}^{0,t,x}\left( \cdot \right) ,0\right) $ is the unique solution
of Eqs. (\ref{2.15}). Therefore, the conclusions in Lemma \ref{l1}, Lemma %
\ref{l2}, and Lemma \ref{l3} also hold when $\varepsilon \rightarrow 0.$

\section{Main Results}

\label{s2}

\subsection{Convergence of Distributions}

In this subsection we first study FBSDEs (\ref{1.1}) with small noise
intensity. To conclude this subsection we introduce the notions of
pseudo-path topology and quasimartingales (Meyer-Zheng \cite{MZ}), adjusted
to our setting. For that, let us introduce the following notations.

Let $T>0$ be a real constant (the terminal time). For any natural number $%
l>0 $, by $D\left( \mathbb{R}^{l}\right) \doteq D\left( \left[ 0,T\right] ;%
\mathbb{R}^{l}\right) $ we denote the Skorohod space of right continuous
with left-hand limits functions $x$ on $\left[ 0,T\right] $ with values in $%
\mathbb{R}^{l}$ such that 
\begin{equation*}
x\left( T-\right) \doteq \lim\limits_{t\uparrow T}x\left( t\right) =x\left(
T\right) ,
\end{equation*}%
and by convention, $x\left( 0-\right) =0.$ Additionally, we introduce the
following metric $\rho $ on $D\left( \mathbb{R}^{l}\right) :$%
\begin{equation*}
\rho \left( x,y\right) =\int_{0}^{T}\left( \left\vert x\left( s\right)
-y\left( s\right) \right\vert \wedge 1\right) \mbox{\rm d}s,\quad x,y\in
D\left( \mathbb{R}^{l}\right) .
\end{equation*}%
The topology induced by this metric is the Meyer-Zheng topology introduced
below on $D\left( \mathbb{R}^{l}\right) .$ Let $\zeta $ be the coordinate
mapping on $D\left( \mathbb{R}^{l}\right) $ defined by 
\begin{equation*}
\zeta _{t}\left( x\right) =x\left( t\right) ,\quad t\in \left[ 0,T\right]
,\quad x\in D\left( \mathbb{R}^{n}\right) ,
\end{equation*}%
and introduce the $\sigma $-algebras of subsets of $D\left( \mathbb{R}%
^{l}\right) ,$%
\begin{equation*}
\mathcal{D}_{s}^{t}\doteq \mathcal{D}_{s}^{t}\left( \mathbb{R}^{l}\right)
\doteq \sigma \left( \left\{ \zeta _{u}:u\in \left[ t,s\right] \right\}
\right) ,\text{ }0\leq t\leq s\leq T,
\end{equation*}%
\begin{equation*}
\mathcal{D\doteq D}\left( \mathbb{R}^{l}\right) \doteq \mathcal{D}%
_{T}^{0}\left( \mathbb{R}^{l}\right) .
\end{equation*}

\noindent \textbf{(Meyer-Zheng topology) }Let $\lambda \left( \text{d}%
t\right) $ the measure $e^{-t}$d$t$ on $\mathbb{R}^{+}.$ Let $w\left(
t\right) $ be a real valued Borel function on $\mathbb{R}_{+}$. The
pseudo-path of $w$ is a probability law on $\left[ 0,\infty \right] \times 
\overline{\mathbb{R}}^{l}:$ the image measure of $\lambda $ under the
mapping $t\rightarrow \left( t,w\left( t\right) \right) .$ We denote by $%
\psi $ the mapping which associates to a path $w$ its pseudo-path: it is
clear that $\psi $ identifies two paths if and only if they are equal a.e.
in Lebesgue's sense. In particular, $\psi $ is $1-1$ on $D\left( \mathbb{R}%
^{l}\right) ,$ and provides an imbedding of $D\left( \mathbb{R}^{l}\right) $
into the compact space $\overline{\mathcal{P}}$ of all probability laws on
the compact space $\left[ 0,\infty \right] \times \overline{\mathbb{R}}^{l}$%
. We give to the induced topology on $D\left( \mathbb{R}^{l}\right) $ the
name of pseudo-path topology or Meyer-Zheng topology. Let us introduce some
intermediate sets between $D\left( \mathbb{R}^{l}\right) ,$ $\overline{%
\mathcal{P}}$ and $\Psi .$ Let $\Psi $ be the set of all pseudo-paths. We
have inclusions 
\begin{equation*}
D\left( \mathbb{R}^{l}\right) \subset \Psi \subset \overline{\mathcal{P}}.
\end{equation*}

The following characterization of the Meyer-Zheng topology is worth noting.

\begin{lemma}
\label{l6}The pseudo-path topology on $\Psi $ is equivalent to the
convergence in measure.
\end{lemma}

Furthermore, it is known that $\Psi $ is a Polish space; and $D\left( 
\mathbb{R}^{l}\right) $ is a Borel set in $\overline{\mathcal{P}}.$

\begin{lemma}
\label{l7}Let $\mathcal{B}\left( D\left( \mathbb{R}^{l}\right) \right) $ be
the $\sigma $-algebra of Borel subsets of $D\left( \mathbb{R}^{l}\right) $
in the Meyer-Zheng topology. Then $\mathcal{B}\left( D\left( \mathbb{R}%
^{l}\right) \right) =\mathcal{D}\left( \mathbb{R}^{l}\right) .$
\end{lemma}

The proof can seen in \cite{BER}. The most important application of the
Meyer-Zheng topology is a tightness result for quasimartingales. We give the
definition here

\begin{definition}
\label{de3}Let $X$ be an $\mathcal{F}$-adapted, cadlag process defined on $%
\left[ 0,T\right] ,$ such that $\mathbb{E}\left\vert X\left( t\right)
\right\vert <\infty $ for all $t\geq 0.$ For any partition $\pi
:0=t_{0}<t_{1}<\cdots <t_{n}\leq T,$ let us define 
\begin{equation}
V_{T}^{\pi }\left( X\right) :=\sum_{0\leq i<n}\mathbb{E}\left\{ \left\vert 
\mathbb{E}\left\{ \left. X\left( t_{i+1}\right) -X\left( t_{i}\right)
\right\vert \mathcal{F}_{t_{i}}\right\} \right\vert \right\} +\mathbb{E}%
\left\vert X\left( t_{n}\right) \right\vert ,  \label{3.1.1}
\end{equation}%
and define the conditional variation of $X$ by $V_{T}\left( X\right)
:=\sup\limits_{\pi }V_{T}^{\pi }\left( X\right) .$ If $V_{T}\left( X\right) $
is finite, then $X$ is called a quasimartingale.
\end{definition}

The following result holds.

\begin{proposition}
\label{pr2}Let $P_{n}$ be a sequence of probability laws on $D\left( \mathbb{%
R}^{l}\right) $ such that under each $P_{n}$ the coordinate process $X\left(
\cdot \right) $ is a quasimartingale with conditional variation $V_{n}\left(
X\left( \cdot \right) \right) $ uniformly bounded in $n.$ Then there exists
a subsequence $\left( P_{n_{k}}\right) $ which converges weakly on $D\left( 
\mathbb{R}^{l}\right) $ to a law $P,$ and $X\left( \cdot \right) $\footnote{%
Note that the quasimartingale in \cite{MZ} is defined on $[0,+\infty ]$.
However, it is fairly easy to check that if $X$ is a quasimartingale on $%
\left[ 0,T\right] $ as is defined above, then the process $X\left( t\right)
=X\left( t\right) \mathbf{I}_{\left[ 0,T\right) }\left( t\right) +X\left(
T\right) \mathbf{I}_{\left[ T,+\infty \right) }\left( t\right) ,$ $t\in $ $%
\left[ 0,+\infty \right] $ is a quasimartingale in the sense of \cite{MZ}.
Furthermore, the conditional variation $V_{T}\left( X\right) $ defined here,
is exactly the same as $V\left( X\right) $ defined in \cite{MZ}. In other
words, our quasimartingale is a \textquotedblleft local\textquotedblright\
version of the one in \cite{MZ}.} is a quasimartingale under $P.$
\end{proposition}

The proof can be seen in \cite{MZ}.

\noindent Now we turn back to Eqs. (\ref{1.1}). Define 
\begin{equation}
\left\{ 
\begin{array}{ccc}
u^{\varepsilon }\left( t,x\right) & = & Y^{\varepsilon ,t,x}\left( t\right) ,
\\ 
v^{\varepsilon }\left( t,x\right) & = & Z^{\varepsilon ,t,x}\left( t\right) ,%
\end{array}%
\right. \text{ for }\left( t,x\right) \in \left[ 0,T\right] \times \mathbb{R}%
.  \label{3.1.2}
\end{equation}%
From the existence and uniqueness of solution for Eq. (\ref{1.1}), we have
the following Markov property%
\begin{equation}
\left\{ 
\begin{array}{ccccc}
u^{\varepsilon }\left( s,X^{\varepsilon ,t,x}\left( s\right) \right) & = & 
Y^{\varepsilon ,s,X^{\varepsilon ,t,x}\left( s\right) }\left( s\right) & = & 
Y^{\varepsilon ,t,x}\left( s\right) \\ 
v^{\varepsilon }\left( s,X^{\varepsilon ,t,x}\left( s\right) \right) & = & 
Z^{\varepsilon ,s,X^{\varepsilon ,t,x}\left( s\right) }\left( s\right) & = & 
Z^{\varepsilon ,t,x}\left( s\right)%
\end{array}%
\right. \text{ a.e., a.e. }s\in \left[ t,T\right] .  \label{3.1.3}
\end{equation}%
In Lemma \ref{l1} we have already seen that the function $u^{\varepsilon }$
is Lipschitz continuous in $x$, uniformly in $t$. With the help of Lemma \ref%
{l3} we also know the continuity property for $u^{\varepsilon }$ in $t$. We
now set 
\begin{equation}
\left\{ 
\begin{array}{crl}
f^{\varepsilon }\left( s,x\right) & = & f\left( s,x,u^{\varepsilon }\left(
s,x\right) \right) , \\ 
\sigma ^{\varepsilon }\left( s,x\right) & = & \sigma \left(
s,x,u^{\varepsilon }\left( s,x\right) \right) , \\ 
g^{\varepsilon }\left( s,x\right) & = & g\left( s,x,u^{\varepsilon }\left(
s,x\right) ,v^{\varepsilon }\left( s,x\right) \right) .%
\end{array}%
\right.  \label{3.1.4}
\end{equation}

We study the properties of $f^{\varepsilon },$ $\sigma ^{\varepsilon }$, $%
f^{0},$ and $\sigma ^{0}$ as follows:

\begin{lemma}
\label{l8}Assume that (A1) and (A2) hold. Then $f^{\varepsilon }$, $\sigma
^{\varepsilon }$, $f^{0},$ and $\sigma ^{0}$ satisfy uniformly Lipschitz
continuous and $f^{\varepsilon }$, $\sigma ^{\varepsilon }$ converge
uniformly $f^{0},$ $\sigma ^{0},$ respectively. Moreover, $f^{\varepsilon }$%
, $\sigma ^{\varepsilon },$ $f^{0}$ and $\sigma ^{0}$ satisfy sublinear
growth.
\end{lemma}

\proof
For any $x,$ $y\in \mathbb{R}^{n},$ $s\in \left[ t,T\right] ,$ we have 
\begin{eqnarray*}
\left\vert f^{\varepsilon }\left( s,x\right) -f^{\varepsilon }\left(
s,y\right) \right\vert &=&\left\vert f\left( s,x,u^{\varepsilon }\left(
s,x\right) \right) -f\left( s,y,u^{\varepsilon }\left( s,y\right) \right)
\right\vert \\
&\leq &C_{1}\left( \left\vert x-y\right\vert +\left\vert u^{\varepsilon
}\left( s,x\right) -u^{\varepsilon }\left( s,y\right) \right\vert \right) \\
&=&C_{1}\left\vert x-y\right\vert +C_{1}\mathcal{C}_{1}\left\vert
x-y\right\vert \\
&\leq &\max \left\{ C_{1},\mathcal{C}_{1}C_{1}\right\} \left\vert
x-y\right\vert ,
\end{eqnarray*}%
by Lemma \ref{l1} and (A1). The same properties for $\sigma ^{\varepsilon }$ 
$f^{0}$, and $f^{0}$ are proved similarly. Next, we show the uniform
convergence, 
\begin{eqnarray*}
\left\vert f^{\varepsilon }\left( s,x\right) -f^{0}\left( s,x\right)
\right\vert &=&\left\vert f\left( s,x,u^{\varepsilon }\left( s,x\right)
\right) -f\left( s,x,u^{0}\left( s,x\right) \right) \right\vert \\
&\leq &C_{1}\left\vert u^{\varepsilon }\left( s,x\right) -u^{0}\left(
s,x\right) \right\vert \leq C_{1}\mathcal{C}_{3}\sqrt{\varepsilon },
\end{eqnarray*}%
by Lemma \ref{l4} and (A1). Once again, the same properties for $\sigma
^{\varepsilon }$ are proved similarly. We are going to prove $b^{\varepsilon
}$ satisfies sublinear growth. We have%
\begin{eqnarray*}
\left\vert f^{\varepsilon }\left( s,x\right) -f\left( s,0,0\right)
\right\vert &=&\left\vert f\left( s,x,u^{\varepsilon }\left( s,x\right)
\right) -f\left( s,0,0\right) \right\vert \\
&\leq &C_{1}\left( \left\vert x\right\vert +\left\vert u^{\varepsilon
}\left( s,x\right) \right\vert \right) \\
&\leq &C_{1}\left( \left\vert x\right\vert +\sqrt{\mathcal{C}_{2}}\left(
1+\left\vert x\right\vert \right) \right) \\
&=&\left( C_{1}+C_{1}\sqrt{\mathcal{C}_{2}}\right) \left\vert x\right\vert
+C_{1}\sqrt{\mathcal{C}_{2}}.
\end{eqnarray*}%
The assumption (A1) yields 
\begin{equation*}
\left\vert f^{\varepsilon }\left( s,x\right) \right\vert \leq \left(
C_{1}+C_{1}\sqrt{\mathcal{C}_{2}}\right) \left\vert x\right\vert +C_{1}\sqrt{%
\mathcal{C}_{2}}+\sup\limits_{0\leq r\leq T}\left\vert f\left( r,0,0\right)
\right\vert .
\end{equation*}%
Once again, the same properties for $\sigma ^{\varepsilon },$ $f^{0},$ $%
\sigma ^{0}$ are proved similarly.%
\endproof%

Then it is easy to check that $\left( X^{\varepsilon ,t,x}\left( \cdot
\right) ,Y^{\varepsilon ,t,x}\left( \cdot \right) ,Z^{\varepsilon
,t,x}\left( \cdot \right) \right) $ solves the following decoupled FBSDEs%
\begin{equation}
\left\{ 
\begin{array}{crl}
X^{\varepsilon ,t,x}\left( s\right) & = & x+\int_{t}^{s}f^{\varepsilon
}\left( r,X^{\varepsilon ,t,x}\left( r\right) \right) \mbox{\rm d}r+\sqrt{%
\varepsilon }\int_{t}^{s}\sigma ^{\varepsilon }\left( r,X^{\varepsilon
,t,x}\left( r\right) \right) \mbox{\rm d}W\left( r\right) , \\ 
Y^{\varepsilon ,t,x}\left( s\right) & = & h\left( X^{\varepsilon ,t,x}\left(
T\right) \right) +\int_{s}^{T}g^{\varepsilon }\left( r,X^{\varepsilon
,t,x}\left( r\right) \right) \mbox{\rm d}r-\int_{s}^{T}Z^{\varepsilon
,t,x}\left( r\right) \mbox{\rm d}W\left( r\right) , \\ 
&  & 0\leq t\leq s\leq T.%
\end{array}%
\right.  \label{3.1.5}
\end{equation}

\noindent From now on, we are concerned on the behavior laws of $\left(
X^{\varepsilon ,t,x},Y^{\varepsilon ,t,x}\right) $ when $\varepsilon
\rightarrow 0$.

\begin{theorem}
\label{t1}Under the assumptions (A1) and (A2), we can conclude the following
results:

i) For all $\delta >0,$ 
\begin{equation}
\lim\limits_{\varepsilon \rightarrow 0}P\left\{ \sup\limits_{t\leq s\leq
T}\left\vert X^{\varepsilon ,t,x}\left( s\right) -\mathcal{X}^{t,x}\left(
s\right) \right\vert >\delta \right\} =0.  \label{3.1.6}
\end{equation}

ii) Let $Q^{\varepsilon }=P\left( \left( Y^{\varepsilon ,t,x}\left( \cdot
\right) \right) ^{-1}\right) $ be the probability measure on $D\left( 
\mathbb{R}^{n}\right) $. Then there exists a subsequence $Q^{\varepsilon
_{n}}$ of $Q^{\varepsilon }$ and a probability law $Q$ on $D\left( \mathbb{R}%
^{n}\right) $ such that $Q^{\varepsilon _{n}}$ converges weakly in the
Meyer-Zheng topology to $Q$ as $n\rightarrow +\infty .$
\end{theorem}

\proof
It follows from the definition of $X^{\varepsilon ,t,x}$ and $\mathcal{X}%
^{t,x}$ that%
\begin{eqnarray}
\sup\limits_{t\leq s\leq T}\left\vert X^{\varepsilon ,t,x}\left( s\right) -%
\mathcal{X}^{t,x}\left( s\right) \right\vert &\leq &\int_{t}^{s}\left\vert
f^{\varepsilon }\left( r,X^{\varepsilon ,t,x}\left( r\right) \right)
-f^{0}\left( r,\mathcal{X}^{t,x}\left( r\right) \right) \right\vert %
\mbox{\rm d}r  \notag \\
&&+\sqrt{\varepsilon }\sup\limits_{t\leq s\leq T}\left\vert
\int_{t}^{s}\sigma \left( r,X^{\varepsilon ,t,x}\left( r\right) \right) %
\mbox{\rm d}W\left( r\right) \right\vert .  \label{3.1.7}
\end{eqnarray}%
From Chebyshev's inequality and the first assertion of the Lemma \ref{l1} we
obtain an estimate of the first term of the right side of (\ref{3.1.7}):%
\begin{eqnarray}
&&\qquad P\left\{ \int_{t}^{s}\left\vert f^{\varepsilon }\left(
r,X^{\varepsilon ,t,x}\left( r\right) \right) -f^{0}\left( r,\mathcal{X}%
^{t,x}\left( r\right) \right) \right\vert \mbox{\rm d}r>\frac{\delta }{2}%
\right\}  \notag \\
&\leq &4\delta ^{-2}\mathbb{E}\left[ \left\vert \int_{t}^{s}\left\vert
f^{\varepsilon }\left( r,X^{\varepsilon ,t,x}\left( r\right) \right)
-f^{0}\left( r,\mathcal{X}^{t,x}\left( r\right) \right) \right\vert %
\mbox{\rm d}r\right\vert ^{2}\right]  \notag \\
&\leq &4\delta ^{-2}T\mathbb{E}\left[ \int_{t}^{s}\left\vert f^{\varepsilon
}\left( r,X^{\varepsilon ,t,x}\left( r\right) \right) -f^{0}\left( r,%
\mathcal{X}^{t,x}\left( r\right) \right) \right\vert ^{2}\mbox{\rm d}r\right]
\notag \\
&=&4\delta ^{-2}T\mathbb{E}\left[ \int_{t}^{s}\left\vert f^{\varepsilon
}\left( r,X^{\varepsilon ,t,x}\left( r\right) \right) -f^{\varepsilon
}\left( r,\mathcal{X}^{t,x}\left( r\right) \right) +f^{\varepsilon }\left( r,%
\mathcal{X}^{t,x}\left( r\right) \right) -f^{0}\left( r,\mathcal{X}%
^{t,x}\left( r\right) \right) \right\vert ^{2}\mbox{\rm d}r\right]  \notag \\
&\leq &8\delta ^{-2}T\mathbb{E}\left[ \int_{t}^{T}\left\vert f^{\varepsilon
}\left( r,X^{\varepsilon ,t,x}\left( r\right) \right) -f^{\varepsilon
}\left( r,\mathcal{X}^{t,x}\left( r\right) \right) \right\vert ^{2}%
\mbox{\rm
d}r+T\mathcal{C}_{3}\varepsilon \right]  \notag \\
&\leq &8\delta ^{-2}T\left( \max \left\{ C_{1},\mathcal{C}_{1}C_{1}\right\}
\right) ^{2}\mathbb{E}\left[ \int_{t}^{T}\left\vert X^{\varepsilon
,t,x}\left( r\right) -\mathcal{X}^{t,x}\left( r\right) \right\vert ^{2}%
\mbox{\rm d}r+T\left( C_{1}\mathcal{C}_{3}\right) ^{2}\varepsilon \right] 
\notag \\
&\leq &8\delta ^{-2}T\left( \max \left\{ C_{1},\mathcal{C}_{1}C_{1}\right\}
\right) ^{2}T\mathbb{E}\left[ \mathcal{C}_{3}\sqrt{\varepsilon }+\left( C_{1}%
\mathcal{C}_{3}\right) ^{2}\varepsilon \right] .  \label{3.1.8}
\end{eqnarray}%
The estimation of the second term in (\ref{3.1.7}) can be accomplished with
the use of the generalized Kolmogorov inequality for stochastic integrals:%
\begin{eqnarray*}
&&\qquad P\left\{ \sqrt{\varepsilon }\sup\limits_{t\leq s\leq T}\left\vert
\int_{t}^{s}\sigma ^{\varepsilon }\left( r,X^{\varepsilon ,t,x}\left(
r\right) \right) \mbox{\rm d}W\left( r\right) \right\vert >\frac{\delta }{2}%
\right\} \\
&\leq &4\delta ^{-2}\varepsilon \mathbb{E}\left[ \int_{t}^{T}\left\vert
\sigma ^{\varepsilon }\left( r,X^{\varepsilon ,t,x}\left( r\right) \right)
\right\vert ^{2}\mbox{\rm d}r\right] \\
&=&4\delta ^{-2}\varepsilon \mathbb{E}\left[ \int_{t}^{T}\left\vert \sigma
\left( r,X^{\varepsilon ,t,x}\left( r\right) ,Y^{\varepsilon ,t,x}\left(
r\right) \right) -\sigma \left( r,0,0\right) +\sigma \left( r,0,0\right)
\right\vert ^{2}\mbox{\rm d}r\right] \\
&\leq &4\delta ^{-2}\varepsilon \mathbb{E}\left[ 4C_{1}^{2}\int_{t}^{T}%
\left( \left\vert X^{\varepsilon ,t,x}\left( r\right) \right\vert
^{2}+\left\vert Y^{\varepsilon ,t,x}\left( r\right) \right\vert
^{2}+2\left\vert \sigma \left( r,0,0\right) \right\vert ^{2}\right) %
\mbox{\rm d}r\right]
\end{eqnarray*}%
\begin{equation}
\leq 4\delta ^{-2}\varepsilon \left[ 8TC_{1}^{2}\mathcal{C}_{2}\left(
1+\left\vert x\right\vert ^{2}\right) +\int_{0}^{T}2\left\vert \sigma \left(
r,0,0\right) \right\vert ^{2}\mbox{\rm d}r\right] ,  \label{3.1.9}
\end{equation}%
where we have used Lemma \ref{l2}. Estimates (\ref{3.1.8})-(\ref{3.1.9})
imply the assertion of the theorem (i).

We are going to prove the second one. First, we establish the connection of
solution of BSDE (\ref{3.1.5}) to quasimartingales. Given a subdivision $\pi
:0=t_{0}<t_{1}<\cdots t_{n}=T$, we get%
\begin{eqnarray}
V_{T}^{\pi }\left( Y^{\varepsilon ,t,x}\right) &\doteq &\mathbb{E}\left\vert
h\left( X^{\varepsilon ,t,x}\left( T\right) \right) \right\vert
+\sum\limits_{k=0}^{n-1}\mathbb{E}\left[ \left\vert \mathbb{E}\left[ \left.
Y^{\varepsilon ,t,x}\left( t_{k+1}\right) -Y^{\varepsilon ,t,x}\left(
t_{k}\right) \right\vert \mathcal{F}_{t}\right] \right\vert \right]  \notag
\\
&=&\mathbb{E}\left\vert h\left( X^{\varepsilon ,t,x}\left( T\right) \right)
\right\vert +\sum\limits_{k=0}^{n-1}\mathbb{E}\left[ \left\vert \mathbb{E}%
\left[ \left. \int_{t_{k}}^{t_{k+1}}g\left( r,X^{\varepsilon ,t,x}\left(
r\right) ,Y^{\varepsilon ,t,x}\left( r\right) ,Z^{\varepsilon ,t,x}\left(
r\right) \right) \mbox{\rm d}r\right\vert \mathcal{F}_{t}\right] \right\vert %
\right]  \notag \\
&\leq &\mathbb{E}\left\vert h\left( X^{\varepsilon ,t,x}\left( T\right)
\right) \right\vert +\mathbb{E}\left[ \int_{0}^{T}g\left( r,X^{\varepsilon
,t,x}\left( r\right) ,Y^{\varepsilon ,t,x}\left( r\right) ,Z^{\varepsilon
,t,x}\left( r\right) \right) \mbox{\rm d}r\right]  \notag \\
&\leq &\mathbb{E}\left\vert h\left( X^{\varepsilon ,t,x}\left( T\right)
\right) -h\left( 0\right) +h\left( 0\right) \right\vert  \notag \\
&&+\mathbb{E}\left[ \int_{t}^{T}\left\vert g^{\varepsilon ,t,x}\left(
r\right) -g\left( r,0,0,0\right) +g\left( r,0,0,0\right) \right\vert %
\mbox{\rm d}r\right]  \notag \\
&\leq &C_{1}\mathbb{E}\left\vert X^{\varepsilon ,t,x}\left( T\right)
\right\vert +C_{1}\mathbb{E}\left[ \int_{t}^{T}\left( \left\vert
X^{\varepsilon ,t,x}\left( r\right) \right\vert +\left\vert Y^{\varepsilon
,t,x}\left( r\right) \right\vert +\left\vert Z^{\varepsilon ,t,x}\left(
r\right) \right\vert \right) \mbox{\rm d}r\right]  \notag \\
&&+\int_{t}^{T}\left\vert g\left( r,0,0,0\right) \right\vert \mbox{\rm d}%
r+h\left( 0\right) .  \label{3.1.10}
\end{eqnarray}%
By Jensen's inequality, it follows from Lemma \ref{l2}%
\begin{equation*}
\mathbb{E}\left\vert X^{\varepsilon ,t,x}\left( T\right) \right\vert \leq 
\sqrt{\mathcal{C}_{2}\left( 1+\left\vert x\right\vert ^{2}\right) },
\end{equation*}%
and 
\begin{equation*}
\begin{array}{l}
\qquad \mathbb{E}\left[ \int_{t}^{T}\left( \left\vert X^{\varepsilon
,t,x}\left( r\right) \right\vert +\left\vert Y^{\varepsilon ,t,x}\left(
r\right) \right\vert +\left\vert Z^{\varepsilon ,t,x}\left( r\right)
\right\vert \right) \mbox{\rm d}rr\right] \\ 
\leq \sqrt{T}\mathbb{E}\left[ \left( \int_{t}^{T}\left\vert X^{\varepsilon
,t,x}\left( r\right) \right\vert ^{2}\mbox{\rm d}r\right) ^{\frac{1}{2}%
}+\left( \int_{t}^{T}\left\vert Y^{\varepsilon ,t,x}\left( r\right)
\right\vert ^{2}\mbox{\rm d}r\right) ^{\frac{1}{2}}+\left(
\int_{t}^{T}\left\vert Z^{\varepsilon ,t,x}\left( r\right) \right\vert ^{2}%
\mbox{\rm d}r\right) ^{\frac{1}{2}}\right] \\ 
\leq \sqrt{T}\left[ \left( \mathbb{E}\int_{t}^{T}\left\vert X^{\varepsilon
,t,x}\left( r\right) \right\vert ^{2}\mbox{\rm d}r\right) ^{\frac{1}{2}%
}+\left( \mathbb{E}\int_{t}^{T}\left\vert Y^{\varepsilon ,t,x}\left(
r\right) \right\vert ^{2}\mbox{\rm d}r\right) ^{\frac{1}{2}}+\left( \mathbb{E%
}\int_{t}^{T}\left\vert Z^{\varepsilon ,t,x}\left( r\right) \right\vert ^{2}%
\mbox{\rm d}r\right) ^{\frac{1}{2}}\right] \\ 
\leq 3\sqrt{T}\left( T\mathcal{C}_{2}\left( 1+\left\vert x\right\vert
^{2}\right) \right) ^{\frac{1}{2}},%
\end{array}%
\end{equation*}%
by Hölder inequality and Jensen's inequality for concave functions. So we
have 
\begin{equation*}
V_{T}^{\pi }\left( Y^{\varepsilon ,t,x}\right) \leq C_{1}\sqrt{\mathcal{C}%
_{2}\left( 1+\left\vert x\right\vert ^{2}\right) }+3\sqrt{T}\left( \mathcal{C%
}_{2}\left( 1+\left\vert x\right\vert ^{2}\right) \right) ^{\frac{1}{2}%
}+h\left( 0\right) +T\sup\limits_{0\leq r\leq T}\left\vert g\left(
r,0,0,0\right) \right\vert .
\end{equation*}%
Hence, by (A1), noting that $V_{T}\left( X\right) \left( Y^{\varepsilon
,t,x}\right) =\sup\limits_{\pi }V_{T}^{\pi }\left( Y^{\varepsilon
,t,x}\right) <+\infty ,$ the result follows.

Now since $D\left( \mathbb{R}^{n}\right) $ is a separable metric space,
there exists a compact metric space $K$ such that $D\left( \mathbb{R}%
^{n}\right) $ is a subset of $K$. Note that $D\left( \mathbb{R}^{n}\right) $
is a Lusin space: for every embedding in a compact metric space $K$, $%
D\left( \mathbb{R}^{n}\right) $ is a Borel set in $K$: $D\left( \mathbb{R}%
^{n}\right) \in \mathcal{B}\left( K\right) $. On the compact metric space $K$
we define 
\begin{equation*}
\widetilde{Q}^{\varepsilon }\left( A\right) =Q^{\varepsilon }\left( A\cap
D\left( \mathbb{R}^{n}\right) \right) ,\quad A\in \mathcal{B}\left( K\right)
.
\end{equation*}%
Clearly, $A\cap D\left( \mathbb{R}^{n}\right) $ belongs to $\mathcal{B}%
\left( D\left( \mathbb{R}^{n}\right) \right) =\mathcal{D}\left( \mathbb{R}%
^{n}\right) $, the last equality being true in view of Lemma \ref{l7}. The
set of probability measures on the compact metric space $K$ is compact for
the weak convergence. Hence, we can choose a subsequence also denoted $%
\left( \varepsilon _{n}\right) _{n\geq 1}$, and a probability measure $%
\widetilde{Q}$ on $K$ such that 
\begin{equation*}
\widetilde{Q}^{\varepsilon _{n}}\overset{w}{\rightarrow }\widetilde{Q}\text{
on }K.
\end{equation*}%
We now show that 
\begin{equation*}
\widetilde{Q}\left( D\left( \mathbb{R}^{n}\right) \right) =1.
\end{equation*}%
We notice that $\widetilde{Q}^{\varepsilon }$ is the distribution of $%
Y^{x,\varepsilon }$ considered as a random variable with values in $\left( K,%
\mathcal{B}\left( K\right) \right) .$ Furthermore, by Proposition \ref{pr2},
we know that, possibly along a subsequence, $\widetilde{Q}^{\varepsilon
_{n}} $ converges weakly to a probability law $Q^{\ast }\in \mathcal{M}%
\left( D\left( \mathbb{R}^{n}\right) \right) .$ The uniqueness of the weak
limit implies that 
\begin{equation}
Q^{\ast }\left( A\right) =\widetilde{Q}\left( A\right) ,\quad \forall A\in 
\mathcal{B}\left( D\left( \mathbb{R}^{n}\right) \right) ,  \label{3.1.11}
\end{equation}%
In particular, 
\begin{equation*}
1=Q^{\ast }\left( D\left( \mathbb{R}^{n}\right) \right) =\widetilde{Q}\left(
D\left( \mathbb{R}^{n}\right) \right) .
\end{equation*}%
The proof is complete. 
\endproof%

\subsection{Almost Sure Convergence}

Considering assumptions (A3)-(A4), we have the following

\begin{theorem}
\label{t2}Under the assumptions (A3) and (A4), supposing further $T\leq C$,
with $C=C(C_{1},T)$ depending on the Lipschitz constant $C_{1}$ and on $T$,
we have

\noindent 1. For each $s,t\in \lbrack 0,T]$, $t\leq s$, the solution of (\ref%
{2.1}), $(X^{\varepsilon ,t,x}\left( s\right) ,Y^{\varepsilon ,t,,x}\left(
s\right) ,Z^{\varepsilon ,t,x}\left( s\right) )_{t\leq s\leq T}$ converges
in $\mathcal{N}_{t}$, when $\varepsilon \rightarrow 0$ to $(X\left( s\right)
,Y\left( s\right) ,0)_{t\leq s\leq T}$, where $(X\left( s\right) ,Y\left(
s\right) )_{t\leq s\leq T}$ solves the coupled system of differential
equations:%
\begin{equation}
\left\{ 
\begin{array}{rcl}
\dot{X}\left( s\right) & = & f(s,X\left( s\right) ,Y\left( s\right) ), \\ 
\dot{Y}\left( s\right) & = & -g(s,X\left( s\right) ,Y\left( s\right)
,0),\,\,t\leq s\leq T, \\ 
X\left( t\right) & = & x,\text{ }Y\left( T\right) =h(X\left( T\right) ).%
\end{array}%
\right.  \label{3.2.1}
\end{equation}%
2. Denoting $u(t,x)=Y^{t,x}\left( t\right) $ the limit in $\varepsilon
\rightarrow 0$ of $Y^{\varepsilon ,t,x}\left( t\right) $, the function $u$
is a viscosity solution of%
\begin{equation}
\begin{cases}
\dfrac{\partial u^{l}}{\partial t}+\sum_{i=1}^{d}f_{i}(t,x,u(t,x))\dfrac{%
\partial u^{l}}{\partial x_{i}}(t,x)+g^{l}(t,x,u(t,x),0)=0, \\ 
u^{\varepsilon }(t,x)=h(x);\,\,x\in \,\mathbb{R}^{d};\,\,t\in \lbrack
0,T];\,\,l=1,...,n.%
\end{cases}
\label{3.2.2}
\end{equation}%
3. The function $u$ is bounded, continuous Lipschitz in $x$ and uniformly
continuous in time.

\noindent 4. Furthermore, if $u\in C_{b}^{1,1}([0,T]\times \mathbb{R}^{n})$,
since (\ref{3.2.1}) has a unique continuous solution, the function $u$ is a
classical solution of (\ref{3.2.2}).
\end{theorem}

\proof
Given $\varepsilon ,\varepsilon _{1}>0$, $x\in \mathbb{R}^{n}$, for $T\leq C$%
, if $(X_{s}^{\varepsilon ,t,x},Y_{s}^{\varepsilon ,t,x},Z_{s}^{t,x})_{t\leq
s\leq T}$ is the unique solution in $\mathcal{N}_{t}$ of%
\begin{equation}
\left\{ 
\begin{array}{rcl}
X^{t,\varepsilon ,x}\left( s\right) & = & x+\int_{t}^{s}f(r,X^{\varepsilon
,t,x}\left( r\right) ,Y^{\varepsilon ,t,x}\left( r\right) )\mbox{\rm d}r \\ 
&  & +\sqrt{\varepsilon }\int_{t}^{s}\sigma (r,X^{\varepsilon ,t,x}\left(
r\right) ,Y^{\varepsilon ,t,x}\left( r\right) )\mbox{\rm d}W(r), \\ 
Y^{\varepsilon ,t,x}\left( s\right) & = & h(X^{\varepsilon ,t,x}\left(
T\right) )+\int_{s}^{T}g(r,X^{\varepsilon ,t,x}\left( r\right)
,Y^{\varepsilon ,t,x}\left( r\right) ,Z^{\varepsilon ,t,x}\left( r\right) )%
\mbox{\rm d}r \\ 
&  & -\int_{s}^{T}Z^{\varepsilon ,t,x}\left( r\right) \mbox{\rm d}W(r),\text{
}x\in \,\mathbb{R}^{d};\,\,\,0\leq t\leq s\leq T,%
\end{array}%
\right.  \label{3.2.3}
\end{equation}%
and $(X_{s}^{\varepsilon _{1},t,x},Y_{s}^{\varepsilon
_{1},t,x},Z_{s}^{\varepsilon _{1},t,x})_{t\leq s\leq T}$ is the unique
solution in $\mathcal{N}_{t}$ of%
\begin{equation}
\left\{ 
\begin{array}{rcl}
X^{\varepsilon _{1},t,x}\left( s\right) & = & x+\int_{t}^{s}f(r,X^{%
\varepsilon _{1},t,x}\left( r\right) ,Y^{\varepsilon _{1},t,x})\mbox{\rm d}r
\\ 
&  & +\sqrt{\varepsilon _{1}}\int_{t}^{s}\sigma (r,X^{\varepsilon
_{1},t,x}\left( r\right) ,Y^{\varepsilon _{1},t,x}\left( r\right) )%
\mbox{\rm
d}W(r), \\ 
Y^{\varepsilon _{1},t,x}\left( s\right) & = & h(X^{\varepsilon
_{1},t,x}\left( T\right) )+\int_{s}^{T}g(r,X^{\varepsilon _{1},t,x}\left(
r\right) ,Y^{\varepsilon _{1},t,x}\left( r\right) ,Z^{\varepsilon
_{1},t,x}\left( r\right) )\mbox{\rm
d}r \\ 
&  & -\int_{s}^{T}Z^{\varepsilon _{1},t,x}\left( r\right) \mbox{\rm d}W(r),%
\text{ }x\in \mathbb{R}^{d};\,\,\,0\leq t\leq s\leq T.%
\end{array}%
\right.  \label{3.2.4}
\end{equation}%
Write, for simplicity of the notations, with $\delta =\varepsilon ,$ $%
\varepsilon _{1}$%
\begin{eqnarray*}
f^{\delta }(s) &=&f(s,X^{\delta ,t,x}\left( s\right) ,Y^{\delta ,t,x}\left(
s\right) ), \\
g^{\delta }(s) &=&g(s,X^{\delta ,t,x}\left( s\right) ,Y^{\delta ,t,x}\left(
s\right) ,Z^{\delta ,t,x}\left( s\right) ), \\
\sigma ^{\delta }(s) &=&\sqrt{\delta }\sigma (s,X^{\delta ,t,x}\left(
s\right) ,Y^{\delta ,t,x}\left( s\right) ), \\
h^{\delta } &=&h(X^{\delta ,t,x}\left( T\right) ).
\end{eqnarray*}%
As in the proof of Theorem 1.2 of \cite{De}, using\textit{\ }Itô's formula:%
\begin{eqnarray}
&&\mathbb{E}\left[ \sup_{t\leq s\leq T}\left\vert X^{\varepsilon }\left(
s\right) -X^{\varepsilon _{1}}\left( s\right) \right\vert ^{2}\right] \leq 
\mathbb{E}\left[ \int_{t}^{T}\left\vert \sigma ^{\varepsilon }(s)-\sigma
^{\varepsilon _{1}}(s)\right\vert ^{2}\mbox{\rm d}s\right]  \notag \\
&&+2\mathbb{E}\,\left[ \sup_{t\leq s\leq T}\int_{t}^{s}\left\langle
X^{\varepsilon }\left( r\right) -X^{\varepsilon _{1}}\left( r\right)
,f^{\varepsilon }(r)-f^{\varepsilon _{1}}(r)\right\rangle \mbox{\rm d}r%
\right]  \notag \\
&&+2\mathbb{E}\,\left[ \sup_{t\leq s\leq T}\int_{t}^{s}\,\left\langle
X^{\varepsilon }\left( r\right) -X^{\varepsilon _{1}}\left( r\right)
,(\sigma ^{\varepsilon }(r)-\sigma ^{\varepsilon _{1}}(r))\mbox{\rm d}%
W(r)\right\rangle \right] .  \label{3.2.5}
\end{eqnarray}%
Using Burkholder-Davis-Gundy's inequalities, there exists $\gamma >0$ such
that:%
\begin{eqnarray}
&&\mathbb{E}\left[ \sup_{t\leq s\leq T}\left\vert X^{\varepsilon }\left(
s\right) -X^{\varepsilon _{1}}\left( s\right) \right\vert ^{2}\right] \leq 
\mathbb{E}\left[ \int_{t}^{T}\left\vert \sigma ^{\varepsilon }(s)-\sigma
^{\varepsilon _{1}}(s)\right\vert ^{2}\mbox{\rm d}s\right]  \notag \\
&&+2\mathbb{E}\,\left[ \sup_{t\leq s\leq T}\int_{t}^{s}\left\langle
X^{\varepsilon }\left( r\right) -X^{\varepsilon _{1}}\left( r\right)
,f^{\varepsilon }(r)-f^{\varepsilon _{1}}(r)\right\rangle \mbox{\rm d}r%
\right]  \notag \\
&&+2\gamma \mathbb{E}\,\left[ \int_{t}^{T}\left\vert X^{\varepsilon }\left(
r\right) -X^{\varepsilon _{1}}\left( r\right) \right\vert ^{2}\left\vert
\sigma ^{\varepsilon }(r)-\sigma ^{\varepsilon _{1}}(r)\right\vert ^{2}%
\mbox{\rm d}r\right] ^{1/2}.  \label{3.2.6}
\end{eqnarray}

\noindent where $\gamma >0$ only depends on the Lipschitz constant $C_{1}$
and $T$. Using the Lipschitz property (A.1), there exists $\gamma >0$,
eventually different, but only dependent on $C_{1},$ $T$ such that:%
\begin{eqnarray}
&&\mathbb{E}\left[ \sup_{t\leq s\leq T}\left\vert X^{\varepsilon }\left(
s\right) -X^{\varepsilon _{1}}\left( s\right) \right\vert ^{2}\right]  \notag
\\
&\leq &\gamma \Bigg \{\mathbb{E}\int_{t}^{T}(\left\vert X^{\varepsilon
}\left( s\right) -X^{\varepsilon _{1}}\left( s\right) \right\vert
^{2}+\left\vert Y^{\varepsilon }\left( s\right) -Y^{\varepsilon _{1}}\left(
s\right) \right\vert ^{2})\mbox{\rm d}s  \notag \\
&&+\mathbb{E}\int_{t}^{T}\left\vert \sigma ^{\varepsilon }(s)-\sigma
^{\varepsilon _{1}}(s)\right\vert ^{2}\mbox{\rm d}s+2\left\vert \sqrt{%
\varepsilon }-\sqrt{\varepsilon _{1}}\right\vert ^{2}  \notag \\
&&+2\sqrt{T}\sup_{t\leq s\leq T}\left\vert X^{\varepsilon }\left( s\right)
-X^{\varepsilon _{1}}\left( s\right) \right\vert ^{2}\Bigg \}.  \label{3.2.7}
\end{eqnarray}%
Then, assuming that $1-2\gamma \sqrt{T}\geq 0,$ we have%
\begin{eqnarray}
&&(1-2\gamma \sqrt{T})\mathbb{E}\left[ \sup_{t\leq s\leq T}\left\vert
X^{\varepsilon }\left( s\right) -X^{\varepsilon _{1}}\left( s\right)
\right\vert ^{2}\right]  \notag \\
&\leq &\gamma \Bigg \{\mathbb{E}\Bigg [\int_{t}^{T}(\left\vert
X^{\varepsilon }\left( s\right) -X^{\varepsilon _{1}}\left( s\right)
\right\vert ^{2}+\left\vert Y^{\varepsilon }\left( s\right) -Y^{\varepsilon
_{1}}\left( s\right) \right\vert ^{2})\mbox{\rm d}s  \notag \\
&&+2\left\vert \sqrt{\varepsilon }-\sqrt{\varepsilon _{1}}\right\vert
^{2}+\int_{t}^{T}\left\vert \sigma ^{\varepsilon }(s)-\sigma ^{\varepsilon
_{1}}(s)\right\vert ^{2}\mbox{\rm d}s\Bigg ]\Bigg \}.  \label{3.2.8}
\end{eqnarray}

\noindent As before, we get 
\begin{eqnarray*}
&&\mathbb{E}\left[ \sup_{t\leq s\leq T}\left\vert Y^{\varepsilon }\left(
s\right) -Y^{\varepsilon _{1}}\left( s\right) \right\vert
^{2}+\int_{t}^{T}\left\vert Z^{\varepsilon }\left( s\right) -Z^{\varepsilon
_{1}}\left( s\right) \right\vert ^{2}\mbox{\rm d}s\right] \\
&\leq &\gamma \Bigg \{\mathbb{E}\left\vert h^{\varepsilon }-h^{\varepsilon
_{1}}\right\vert ^{2}+\mathbb{E}\Big [\sup_{t\leq s\leq
T}\int_{t}^{T}\left\langle Y^{\varepsilon }\left( s\right) -Y^{\varepsilon
_{1}}\left( s\right) ,g^{\varepsilon }(s)-g^{\varepsilon
_{1}}(s)\right\rangle \mbox{\rm d}s\Big ]\Bigg \}.
\end{eqnarray*}

\noindent Using (\ref{3.2.7}) and the assumption (A.3) modifying $\gamma $
eventually:%
\begin{eqnarray*}
&&\mathbb{E}\left[ \sup_{t\leq s\leq T}\left\vert X^{\varepsilon }\left(
s\right) -X^{\varepsilon _{1}}\left( s\right) \right\vert ^{2}+\sup_{t\leq
s\leq T}\left\vert Y^{\varepsilon }\left( s\right) -Y^{\varepsilon
_{1}}\left( s\right) \right\vert ^{2}+\int_{t}^{T}\left\vert Z^{\varepsilon
}\left( s\right) -Z^{\varepsilon _{1}}\left( s\right) \right\vert ^{2}%
\mbox{\rm d}s\right] \\
&\leq &\gamma \Bigg \{\mathbb{E}\left\vert h^{\varepsilon }-h^{\varepsilon
_{1}}\right\vert ^{2}+\mathbb{E}\int_{t}^{T}\left\vert X^{\varepsilon
}\left( s\right) -X^{\varepsilon _{1}}\left( s\right) \right\vert ^{2}%
\mbox{\rm d}s+\mathbb{E}\int_{t}^{T}\left\vert Y^{\varepsilon }\left(
s\right) -Y^{\varepsilon _{1}}\left( s\right) \right\vert ^{2}\mbox{\rm d}s
\\
&&+\mathbb{E}\Big (\int_{t}^{T}(\left\vert f^{\varepsilon
}(s)-f^{\varepsilon _{1}}(s)\right\vert +\left\vert g^{\varepsilon
}(s)-g^{\varepsilon _{1}}(s)\right\vert )\mbox{\rm d}s\Big )^{2} \\
&&+\mathbb{E}\int_{t}^{T}\left\vert \sigma ^{\varepsilon }(s)-\sigma
^{\varepsilon _{1}}(s)\right\vert ^{2}\mbox{\rm d}s+\left\vert \sqrt{%
\varepsilon }-\sqrt{\varepsilon _{1}}\right\vert ^{2} \\
&&+\mathbb{E}\int_{t}^{T}\left\vert Z^{\varepsilon }\left( s\right)
-Z^{\varepsilon _{1}}\left( s\right) \right\vert (\left\vert Y^{\varepsilon
}\left( s\right) -Y^{\varepsilon _{1}}\left( s\right) \right\vert
+\left\vert X^{\varepsilon }\left( s\right) -X^{\varepsilon _{1}}\left(
s\right) \right\vert )\mbox{\rm d}s\Bigg \}.
\end{eqnarray*}

\noindent So there exist $C^{\ast }\leq C$ and $\Gamma $, only depending on $%
C_{1}$ and $\gamma _{1}$, such that, for $T-t\leq C^{\ast }$%
\begin{eqnarray}
&&\mathbb{E}\left[ \sup_{t\leq s\leq T}\left\vert X^{\varepsilon }\left(
s\right) -X^{\varepsilon _{1}}\left( s\right) \right\vert ^{2}+\sup_{t\leq
s\leq T}\left\vert Y^{\varepsilon }\left( s\right) -Y^{\varepsilon
_{1}}\left( s\right) \right\vert ^{2}+\int_{t}^{T}\left\vert Z^{\varepsilon
}\left( s\right) -Z^{\varepsilon _{1}}\left( s\right) \right\vert ^{2}%
\mbox{\rm d}s\right]  \notag \\
&\leq &\gamma \Bigg \{\mathbb{E}\left\vert h^{\varepsilon }-h^{\varepsilon
_{1}}\right\vert ^{2}+\mathbb{E}\int_{t}^{T}\left\vert \sigma ^{\varepsilon
}(s)-\sigma ^{\varepsilon _{1}}(s)\right\vert ^{2}\mbox{\rm d}s+\left\vert 
\sqrt{\varepsilon }-\sqrt{\varepsilon _{1}}\right\vert ^{2}  \notag \\
&&+\mathbb{E}\big (\int_{t}^{T}\left\vert Y^{\varepsilon }\left( s\right)
-Y^{\varepsilon _{1}}\left( s\right) \right\vert +\left\vert g^{\varepsilon
}(s)-g^{\varepsilon _{1}}(s)\right\vert \mbox{\rm d}s\big )^{2}\Bigg \}.
\label{3.2.9}
\end{eqnarray}

\noindent We can repeat the argument in $[T-2C^{\ast };T-C^{\ast }]$ and
recurrently we get this important a-priori estimate (\ref{3.2.9}) valid for
all $[t,T]$, with some (possibly different) constant $\gamma >0$ only
depending on $C_{1}$.

Now, from (A.3), the following estimate holds:

\begin{equation}
\mathbb{E}\left[ \left\vert h^{\varepsilon }-h^{\varepsilon _{1}}\right\vert
^{2}\right] \leq L^{2}\mathbb{E}\left[ \sup\limits_{t\leq s\leq T}\left\vert
X^{\varepsilon }\left( s\right) -X^{\varepsilon _{1}}\left( s\right)
\right\vert ^{2}\right] .  \label{3.2.10}
\end{equation}%
One can also derive the estimate 
\begin{eqnarray}
&&\mathbb{E}\left[ \int_{t}^{T}\Big (\left\vert \sqrt{\varepsilon }\sigma
^{\varepsilon }\left( s\right) -\sqrt{\varepsilon _{1}}\sigma ^{\varepsilon
_{1}}(s)\right\vert ^{2}\Big )\mbox{\rm d}s\right]  \notag \\
&\leq &c_{1}\left\vert \sqrt{\varepsilon }-\sqrt{\varepsilon _{1}}%
\right\vert ^{2}+c_{2}\mathbb{E}\left[ \int_{t}^{T}\big (\left\vert
X^{\varepsilon }\left( s\right) -X^{\varepsilon _{1}}\left( s\right)
\right\vert ^{2}+\left\vert Y^{\varepsilon }\left( s\right) -Y^{\varepsilon
_{1}}\left( s\right) \right\vert ^{2}\big )\mbox{\rm d}s\right] .
\label{3.2.11}
\end{eqnarray}%
where $c_{1},$ $c_{2}>0$ only depending on $T$, $L$ and on the bound of $%
\left\vert \sigma \right\vert $.

Furthermore, using Jensen's inequality and the Lipschitz property of $g$,
for some $c_{3}$ modified eventually along the various steps,%
\begin{eqnarray}
&&\mathbb{E}\left[ \int_{t}^{T}\left\vert Y^{\varepsilon }\left( s\right)
-Y^{\varepsilon _{1}}\left( s\right) \right\vert +\left\vert g^{\varepsilon
}(s)-g^{\varepsilon _{1}}(s)\right\vert \mbox{\rm d}s\right] ^{2}  \notag \\
&\leq &\mathbb{E}\left[ \int_{t}^{T}\left\vert Y^{\varepsilon }\left(
s\right) -Y^{\varepsilon _{1}}\left( s\right) \right\vert ^{2}+\left\vert
g^{\varepsilon }\left( s\right) -g^{\varepsilon _{1}}\left( s\right)
\right\vert ^{2}+2\left\vert Y^{\varepsilon }\left( s\right) -Y^{\varepsilon
_{1}}\left( s\right) \right\vert \left\vert g^{\varepsilon }\left( s\right)
-g^{\varepsilon _{1}}\left( s\right) \right\vert \mbox{\rm d}s\right]  \notag
\\
&\leq &c_{3}\mathbb{E}\Bigg [\int_{t}^{T}\left\vert X^{\varepsilon }\left(
s\right) -X^{\varepsilon _{1}}\left( s\right) \right\vert ^{2}+\left\vert
Y^{\varepsilon }\left( s\right) -Y^{\varepsilon _{1}}\left( s\right)
\right\vert ^{2}+\left\vert Z^{\varepsilon }\left( s\right) -Z^{\varepsilon
_{1}}\left( s\right) \right\vert ^{2}\mbox{\rm d}s\Bigg ]  \notag \\
&\leq &c_{3}\Bigg \{\mathbb{E}\left[ \sup_{t\leq s\leq T}\left\vert
X^{\varepsilon }\left( s\right) -X^{\varepsilon _{1}}\left( s\right)
\right\vert ^{2}\right] +\mathbb{E}\left[ \sup_{t\leq s\leq T}\left\vert
Y^{\varepsilon }\left( s\right) -Y^{\varepsilon _{1}}\left( s\right)
\right\vert ^{2}\right]  \notag \\
&&+\mathbb{E}\left[ \sup_{t\leq s\leq T}\left\vert Z^{\varepsilon }\left(
s\right) -Z^{\varepsilon _{1}}\left( s\right) \right\vert ^{2}\mbox{\rm d}s%
\right] \Bigg \}.  \label{3.2.12}
\end{eqnarray}

\noindent Using (\ref{2.10}) and the boundedness of $\left\vert \sigma
\right\vert $, we have $\mathbb{E}\left[ \int_{t}^{T}\left\vert
Z^{\varepsilon }\left( s\right) \right\vert ^{2}\mbox{\rm d}s\right]
\rightarrow 0$ as $\varepsilon \rightarrow 0$ and

\begin{equation}
\lim_{\left\vert \varepsilon -\varepsilon _{1}\right\vert \rightarrow 0}%
\mathbb{E}\left[ \int_{t}^{T}\left\vert Z^{\varepsilon }\left( s\right)
-Z^{\varepsilon _{1}}\left( s\right) \right\vert ^{2}\mbox{\rm d}s\right] =0.
\label{3.2.13}
\end{equation}%
Furthermore, by Burkholder-Davis-Gundy's inequalities, there exists a
constant $\gamma >0$, eventually new such that:%
\begin{eqnarray}
&&\mathbb{E}\left[ \sup_{t\leq s\leq T}\left\vert X^{\varepsilon }\left(
s\right) -X^{\varepsilon _{1}}\left( s\right) \right\vert ^{2}+\sup_{t\leq
s\leq T}\left\vert Y^{\varepsilon }\left( s\right) -Y^{\varepsilon
_{1}}\left( s\right) \right\vert ^{2}+\int_{t}^{T}\left\vert Z^{\varepsilon
}\left( s\right) -Z^{\varepsilon _{1}}\left( s\right) \right\vert ^{2}%
\mbox{\rm d}s\right]  \notag \\
&\leq &\gamma \mathbb{E}\left[ \int_{t}^{T}\left\vert X^{\varepsilon }\left(
s\right) -X^{\varepsilon _{1}}\left( s\right) \right\vert ^{2}+\left\vert
Y^{\varepsilon }\left( s\right) -Y^{\varepsilon _{1}}\left( s\right)
\right\vert ^{2}+\left\vert Z^{\varepsilon }\left( s\right) -Z^{\varepsilon
_{1}}\left( s\right) \right\vert ^{2}\mbox{\rm d}s\right]  \notag \\
&\leq &\gamma \mathbb{E}\Bigg [\int_{t}^{T}\sup_{t\leq r\leq s}\left\vert
X^{\varepsilon }\left( r\right) -X^{\varepsilon _{1}}\left( r\right)
\right\vert ^{2}+\sup_{t\leq r\leq s}\left\vert Y^{\varepsilon }\left(
r\right) -Y^{\varepsilon _{1}}\left( r\right) \right\vert ^{2}  \notag \\
&&+\sup_{t\leq r\leq s}\left\vert Z^{\varepsilon }\left( r\right)
-Z^{\varepsilon _{1}}\left( r\right) \right\vert ^{2}\mbox{\rm d}s\Bigg ].
\label{3.2.14}
\end{eqnarray}

\noindent Moreover, for some $\gamma _{1},$ $\gamma _{2}>0$, using (\ref%
{3.2.10})-(\ref{3.2.14}), we have

\begin{eqnarray}
&&\mathbb{E}\left[ \sup_{t\leq s\leq T}\left\vert X^{\varepsilon }\left(
s\right) -X^{\varepsilon _{1}}\left( s\right) \right\vert ^{2}+\sup_{t\leq
s\leq T}\left\vert Y^{\varepsilon }\left( s\right) -Y^{\varepsilon
_{1}}\left( s\right) \right\vert ^{2}+\int_{t}^{T}\left\vert Z^{\varepsilon
}\left( s\right) -Z^{\varepsilon _{1}}\left( s\right) \right\vert ^{2}%
\mbox{\rm d}s\right]  \notag \\
&\leq &\gamma _{1}\left\vert \sqrt{\varepsilon }-\sqrt{\varepsilon _{1}}%
\right\vert ^{2}+\gamma _{2}\mathbb{E}\Bigg [\int_{t}^{T}\sup_{t\leq r\leq
s}\left\vert X^{\varepsilon }\left( r\right) -X^{\varepsilon _{1}}\left(
r\right) \right\vert ^{2}\mbox{\rm d}r  \notag \\
&&+\sup_{t\leq r\leq s}\left\vert Y^{\varepsilon }\left( r\right)
-Y^{\varepsilon _{1}}\left( r\right) \right\vert ^{2}+\sup_{t\leq r\leq
s}\left\vert Z^{\varepsilon }\left( r\right) -Z^{\varepsilon _{1}}\left(
r\right) \right\vert ^{2}\mbox{\rm d}s\Bigg ].  \label{3.2.15}
\end{eqnarray}%
Using Gronwall's inequality, we get 
\begin{eqnarray}
&&\mathbb{E}\left[ \int_{t}^{T}\sup_{t\leq r\leq s}\left\vert X^{\varepsilon
}\left( r\right) -X^{\varepsilon _{1}}\left( r\right) \right\vert
^{2}+\sup_{t\leq r\leq s}\left\vert Y^{\varepsilon }\left( s\right)
-Y^{\varepsilon _{1}}\left( s\right) \right\vert ^{2}\mbox{\rm d}r\right] 
\notag \\
&\leq &C_{2}\left\vert \sqrt{\varepsilon }-\sqrt{\varepsilon _{1}}%
\right\vert ^{2}\mathbb{E}\left[ \int_{t}^{T}\sup_{t\leq s\leq T}\left\vert
Z^{\varepsilon }\left( s\right) -Z^{\varepsilon _{1}}\left( s\right)
\right\vert ^{2}\mbox{\rm d}s\right] \rightarrow 0.  \label{3.2.16}
\end{eqnarray}%
for some $C_{2}>0$ only depending on $C_{1}$ (independent of $\varepsilon $%
), since $\mathbb{E}\left[ \sup_{t\leq s\leq T}\left\vert Z^{\varepsilon
}\left( s\right) -Z^{\varepsilon _{1}}\left( s\right) \right\vert ^{2}\right]
$ is bounded, by the results (\ref{2.9}) and (\ref{2.10}) .

We conclude, by the previous estimates, that the pair $(X^{t,\varepsilon
}\left( s\right) ,Y^{t,\varepsilon }\left( s\right) )_{t\leq s\leq T}$ form
a Cauchy sequence and therefore converges in $\mathcal{S}^{2}\left( t,T;%
\mathbb{R}^{d}\right) \times \mathcal{S}^{2}\left( t,T;\mathbb{R}^{k}\right) 
$.

Denote $(X\left( s\right) ,Y\left( s\right) )_{t\leq s\leq T}$ its limit. $%
(X\left( s\right) ,Y\left( s\right) ,0)_{t\leq s\leq T}$ is the limit in $%
\mathcal{N}_{t}$ of 
\begin{equation*}
(X^{t,\varepsilon }\left( s\right) ,Y^{t,\varepsilon }\left( s\right)
,Z^{t,\varepsilon ,x}\left( s\right) )_{t\leq s\leq T}
\end{equation*}
when $\varepsilon \rightarrow 0$.

If we consider the forward equation in (\ref{1.1}) and if we take the limit
pointwise when $\varepsilon \rightarrow 0$, we have

\begin{equation}
X^{t}\left( s\right) =x+\int_{t}^{s}f(r,X\left( r\right) ,Y\left( r\right) )%
\mbox{\rm d}r,  \label{3.2.17}
\end{equation}%
where we have used the boundedness of $\sigma $ and the continuity of $f$.

Similarly, we can take the limit on the backward equation when $\varepsilon
\rightarrow 0$. Using the continuity of the functions $h,$ $g$ and%
\begin{equation*}
\mathbb{E}\left\vert \int_{s}^{T}Z^{t,\varepsilon }\left( r\right) 
\mbox{\rm
d}W(r)\right\vert ^{2}=\mathbb{E}\left[ \int_{s}^{T}\left\vert
Z^{t,\varepsilon }\left( r\right) \right\vert ^{2}\mbox{\rm d}r\right]
\rightarrow 0,\text{ as }\varepsilon \rightarrow 0,
\end{equation*}

\noindent which implies $\int_{s}^{T}Z^{t,\varepsilon }\left( r\right) %
\mbox{\rm d}W(r)\rightarrow 0$, $P$-a.s., we have 
\begin{equation}
Y_{s}^{t}=h(X\left( T\right) )+\int_{s}^{T}g(r,X\left( r\right) ,Y\left(
r\right) ,0)\mbox{\rm d}r.  \label{3.2.18}
\end{equation}%
In conclusion, $(X\left( s\right) ,Y\left( s\right) )_{t\leq s\leq T}$
solves the following deterministic problem of ordinary (coupled)
differential equations, almost surely,%
\begin{equation}
\left\{ 
\begin{array}{rcl}
\dot{X}\left( s\right) & = & f(s,X\left( s\right) ,Y\left( s\right) ), \\ 
\dot{Y}\left( s\right) & = & -g(s,X\left( s\right) ,Y\left( s\right)
,0),\,\,t\leq s\leq T, \\ 
X\left( t\right) & = & x,Y\left( T\right) =h(X\left( T\right) ).%
\end{array}%
\right.  \label{3.2.19}
\end{equation}

Given $t,$ $t^{\prime }\in \lbrack 0,T]$, $x,$ $x^{\prime }\in \mathbb{R}%
^{n} $, consider $(X_{s}^{\varepsilon ,t,x},Y_{s}^{\varepsilon
,t,x},Z_{s}^{\varepsilon ,t,x})_{t\leq s\leq T}$ the unique solution in $%
\mathcal{N}_{t}$%
\begin{equation}
\left\{ 
\begin{array}{rcl}
X^{\varepsilon ,t,x}\left( s\right) & = & x+\int_{t}^{s}f(r,X^{\varepsilon
,t,x}\left( r\right) ,Y^{\varepsilon ,t,x}\left( r\right) )\mbox{\rm d}r \\ 
&  & +\sqrt{\varepsilon }\int_{t}^{s}\sigma (r,X^{\varepsilon ,t,x}\left(
r\right) ,Y^{t,\varepsilon }\left( r\right) )\mbox{\rm d}W(r), \\ 
Y^{\varepsilon ,t,x}\left( s\right) & = & h(X^{\varepsilon ,t,x}\left(
T\right) )+\int_{s}^{T}g(r,X^{\varepsilon ,t,x}\left( r\right)
,Y^{\varepsilon ,t,x}\left( r\right) ,Z^{\varepsilon ,t,x}\left( r\right) )%
\mbox{\rm d}r \\ 
&  & -\int_{s}^{T}Z^{\varepsilon ,t,x}\left( r\right) \mbox{\rm d}W(r),\text{
}x\in \mathbb{R}^{n}\text{, }t\leq s\leq T,%
\end{array}%
\right.  \label{3.2.20}
\end{equation}%
\noindent extended to the whole interval $[0,T]$, putting 
\begin{equation}
\forall \,\,0\leq s\leq t\,\,\,\,X^{\varepsilon ,t,x}\left( s\right)
=x;\,\,Y^{\varepsilon ,t,x}\left( s\right) =Y^{\varepsilon ,t,x}\left(
t\right) ;\,\,Z^{\varepsilon ,t,x}\left( s\right) =0,  \label{3.2.21}
\end{equation}%
and consider 
\begin{equation*}
(X^{\varepsilon ,t^{\prime },x^{\prime }}\left( s\right) ,Y^{\varepsilon
,t^{\prime },x^{\prime }}\left( s\right) ,Z^{\varepsilon ,t^{\prime
},x^{\prime }}\left( s\right) )_{t^{\prime }\leq s\leq T}
\end{equation*}%
the unique solution in $\mathcal{N}_{t^{\prime }}$%
\begin{equation}
\left\{ 
\begin{array}{rcl}
X^{\varepsilon ,t^{\prime },x^{\prime }}\left( s\right) & = & x^{\prime
}+\int_{t^{\prime }}^{s}f(r,X^{\varepsilon ,t^{\prime },x^{\prime }}\left(
r\right) ,Y^{\varepsilon ,t^{\prime },x^{\prime }}\left( r\right) )%
\mbox{\rm
d}r \\ 
&  & +\sqrt{\varepsilon }\int_{t^{\prime }}^{s}\sigma (r,X^{\varepsilon
,t^{\prime },x^{\prime }}\left( r\right) ,Y^{\varepsilon ,t^{\prime
},x^{\prime }}\left( r\right) )\mbox{\rm d}W(r), \\ 
Y^{\varepsilon ,t^{\prime },x^{\prime }}\left( s\right) & = & h\left(
X^{\varepsilon ,t^{\prime },x^{\prime }}\left( T\right) \right)
+\int_{s}^{T}g(r,X^{\varepsilon ,t^{\prime },x^{\prime }}\left( r\right)
,Y^{\varepsilon ,t^{\prime },x^{\prime }}\left( r\right) ,Z^{\varepsilon
,t^{\prime },x^{\prime }}\left( r\right) )\mbox{\rm d}r \\ 
&  & -\int_{s}^{T}Z^{\varepsilon ,t^{\prime },x^{\prime }}\left( r\right) %
\mbox{\rm d}W\left( r\right) ,\text{ }x^{\prime }\in \mathbb{R}^{n},\text{ }%
t^{\prime }\leq s\leq T,%
\end{array}%
\right.  \label{3.2.22}
\end{equation}%
extended to the whole interval $[0,T]$, putting 
\begin{equation}
\forall 0\leq s\leq t^{\prime }\,\,\,\,X^{\varepsilon ,t^{\prime },x^{\prime
}}\left( s\right) =x;\,\,Y^{\varepsilon ,t^{\prime },x^{\prime }}\left(
s\right) =Y^{\varepsilon ,t^{\prime },x^{\prime }}\left( t^{\prime }\right)
;\,\,Z^{\varepsilon ,t^{\prime },x^{\prime }}\left( s\right) =0.
\label{3.2.23}
\end{equation}%
Using the estimate (\ref{3.2.9}), such as in the Corollary 1.4 of \cite{De},
we are lead to%
\begin{eqnarray}
&&\mathbb{E}\bigg[\sup_{0\leq s\leq T}\left\vert X^{\varepsilon t,x}\left(
s\right) -X^{\varepsilon ,t^{\prime },x^{\prime }}\left( s\right)
\right\vert ^{2}+\sup_{0\leq s\leq T}\left\vert Y^{\varepsilon t,x}\left(
s\right) -Y^{\varepsilon ,t^{\prime },x^{\prime }}\left( s\right)
\right\vert ^{2}  \notag \\
&&+\int_{0}^{T}\left\vert Z^{\varepsilon ,t,x}\left( s\right)
-Z^{\varepsilon ,t^{\prime },x^{\prime }}\left( s\right) \right\vert ^{2}%
\bigg]  \notag \\
&\leq &\alpha \left\vert x-x^{\prime }\right\vert ^{2}+\beta (1+\left\vert
x\right\vert ^{2})\left\vert t-t^{\prime }\right\vert ^{2},  \label{3.2.24}
\end{eqnarray}%
where $\alpha ,$ $\beta >0$ are constants only depending on $C_{1},$ $%
\Lambda $.

Finally,%
\begin{equation}
\left\vert u^{\varepsilon }\left( t,x\right) -u^{\varepsilon }\left(
t^{\prime },x^{\prime }\right) \right\vert ^{2}\leq \alpha \left\vert
x-x^{\prime }\right\vert ^{2}+\beta (1+\left\vert x\right\vert
^{2})\left\vert t-t^{\prime }\right\vert ^{2},  \label{3.2.25}
\end{equation}

\noindent which proves that $(u^{\varepsilon })_{\varepsilon >0}$ is a
family of equicontinuous maps on every compact set of $[0,T]\times \mathbb{R}%
^{n}$. We can apply Arzela's Ascoli Theorem and conclude that the
convergence of $u^{\varepsilon }$ to $u$, where $u(t,x)=Y_{t}^{t,x}$, is
uniform in $[0,T]\times K$, for every $K$ compact subset of $\mathbb{R}^{d}$.

Taking the limit in $\varepsilon \rightarrow 0$ in (\ref{3.2.25}), we get%
\begin{equation}
\left\vert u\left( t^{\prime },x^{\prime }\right) -u\left( t,x\right)
\right\vert ^{2}\leq \alpha \left\vert x-x^{\prime }\right\vert ^{2}+\beta
(1+\left\vert x\right\vert ^{2})\left\vert t-t^{\prime }\right\vert ^{2},
\label{3.2.26}
\end{equation}%
for all $\ (t,x),(t^{\prime },x^{\prime })\in \lbrack 0,T]\times \mathbb{R}%
^{n}$, which proves that the function $u$, limit in $\varepsilon \rightarrow
0$, is Lipschitz continuous in $x$ and uniformly continuous in $t$. The
boundedness of $u$ is given by (\ref{2.7}), since this bound is uniform in $%
\varepsilon $. Using Theorem 5.1 of \cite{P1} we deduce that $u^{\varepsilon
}$ is a viscosity solution in $[0,T]\times \mathbb{R}^{d}$ of (\ref{2.2}).

Since the coefficients of the quasilinear parabolic system are Lipschitz
continuous, by the property of the compact uniform convergence of viscosity
solutions for quasilinear parabolic equations (see \cite{CIL} for details),
we conclude that $u$ is a viscosity solution of (\ref{3.2.2}).

Moreover, let $v:[0,T]\times \mathbb{R}^{n}\rightarrow \mathbb{R}^{n}$ be a $%
C_{b}^{1,1}([0,T],\mathbb{R}^{n})$ solution, continuous Lipschitz in $x$ and
uniformly continuous in $t$ for (\ref{3.2.2}). Fixing $(t,x)\in \lbrack
0,T]\times \mathbb{R}^{n}$, we take the following function:%
\begin{eqnarray*}
\psi &:&[t,T]\rightarrow \mathbb{R}^{n}, \\
\psi (s) &:&=v(s,X^{t,x}\left( s\right) ).
\end{eqnarray*}%
Computing its time derivative: 
\begin{eqnarray*}
\dfrac{\text{d}\psi }{\text{d}s}(s) &=&\dfrac{\partial v}{\partial s}%
(s,X^{t,x}\left( s\right) )+\sum_{i=1}^{d}\dfrac{\partial v}{\partial x_{i}}%
(s,X^{t,x}\left( s\right) )\dfrac{\partial (X^{t,x}\left( s\right) )}{%
\partial t} \\
&=&\dfrac{\partial v}{\partial s}(s,X^{t,x}\left( s\right) )+\sum_{i=1}^{d}%
\dfrac{\partial v}{\partial x_{i}}(s,X^{t,x}\left( s\right)
)f(s,X^{t,x}\left( s\right) ,Y^{t,x}\left( s\right) ) \\
&=&-g(s,X^{t,x}\left( s\right) ,v(s,X^{t,x}\left( s\right) ),0), \\
\psi (T) &=&v(T,X^{t,x}\left( T\right) )=h(x).
\end{eqnarray*}

As a consequence, $v(t,x)=v(t,X^{t,x}\left( t\right) )=u(t,x)$, under the
hypothesis of uniqueness of solution for the system of ordinary differential
equations (\ref{3.2.1}). So, under these hypothesis, we have a uniqueness
property of solution for (\ref{3.2.2}) in the class of $C_{b}^{1,1}([0,T]%
\times \mathbb{R}^{n})$, which are Lipschitz continuous in $x$ and uniformly
continuous in $t$. 
\endproof%

\begin{remark}
\label{re2}\textsl{The result presented above (Theorem \ref{t2}) is stated
under the assumptions (A.3) and (A.4), which ensure existence and uniqueness
of solution of (\ref{3.2.3}) for a local time }$T\leq C$\textsl{, where }$C$%
\textsl{\ is a certain constant only depending on the Lipschitz constant }$%
C_{1}$\textsl{. Under these assumptions, our results and the existence and
uniqueness of solution for the FBSDEs (\cite{De, MPY}) do not depend on
results for PDE's but only on probabilistic arguments. However, it is
possible to extend Theorem \ref{t2} to global time. In order to do it and to
prove the important properties in \cite{De}, results of deterministic
quasilinear parabolic partial differential equations \cite{LS} are used. If
we maintain the assumption of smoothness on the coefficients of the FBSDEs
(A.3) and (A.4), the\ Four Step Scheme Methodology ensures existence and
uniqueness of solution for (\ref{1.1}), and if we remove this requirement of
smoothness on the coefficients of the system (\ref{1.1}), a regularization
argument in \cite{De} is used in order to prove it. In this setting, we can
also conclude the claims (2) and (3) of Theorem \ref{t2}.}
\end{remark}

\section{Large Deviation Principle}

\label{s3}

In this section, we study the Freidlin-Wentzell's large deviation principle
for the laws of the family of processes 
\begin{equation*}
\left\{ \left( X^{\varepsilon ,t,x}\left( \cdot \right) ,Y^{\varepsilon
,t,x}\left( \cdot \right) \right) \right\} _{\varepsilon \in \left( 0,1 
\right] }\in \mathcal{S}^{2}\left( t,T;\mathbb{R}^{n}\right) \times \mathcal{%
S}^{2}\left( t,T;\mathbb{R}^{n}\right)
\end{equation*}%
as $\varepsilon \rightarrow 0.$

To begin, let us recall the following definitions mainly from \cite{DZ}.

\begin{definition}
\label{de4}If $E$ is a complete separable metric space, then a function $%
\mathcal{I}$ defined on $E$ is called a rate function if it has the
following properties: 
\begin{equation}
\left\{ 
\begin{array}{l}
\text{(a)\qquad }\mathcal{I}:E\rightarrow \left[ 0,+\infty \right] ,\text{ }%
\mathcal{I}\text{ is lower semicontinuous.} \\ 
\text{(b)\qquad If }0\leq a\leq \infty ,\text{ then }C_{\mathcal{I}}\left(
a\right) =\left\{ x\in E:I\left( x\right) \leq a\right\} \text{ is compact.}%
\end{array}%
\right.  \tag{4.1}
\end{equation}
\end{definition}

\begin{definition}[\textbf{Large deviation principle}]
If $E$ is a complete separable metric space, $\mathcal{B}$ is the Borel $%
\sigma $-field on $E,$ $\left\{ \mu _{\varepsilon }:\varepsilon >0\right\} $
is a family of probability measures on $\left( E,\mathcal{B}\right) ,$ and $%
\mathcal{I}$ is a function defined on $E$ and satisfying (4.1), then we say
that $\left\{ \mu _{\varepsilon }\right\} _{\varepsilon >0}$ satisfies a
large deviation principle with rate $\mathcal{I}$ if: 
\begin{equation}
\left\{ 
\begin{array}{l}
\text{(a)\qquad For every open subset }A\text{ of }E\text{,} \\ 
\qquad \qquad \liminf_{\varepsilon \rightarrow 0}\varepsilon \log \mu
_{\varepsilon }\left( A\right) \geq -\inf\limits_{g\in A}\mathcal{I}\left(
g\right) . \\ 
\text{(b)\qquad For every closed subset }A\text{ of }E, \\ 
\qquad \qquad \limsup_{\varepsilon \rightarrow 0}\varepsilon \log \mu
_{\varepsilon }\left( A\right) \leq -\inf\limits_{g\in A}\mathcal{I}\left(
g\right) ,%
\end{array}%
\right.  \tag{4.2}
\end{equation}%
here the infimum over the empty set is defined to be +$\infty .$
\end{definition}

For our purpose recall the following results from \cite{A,BC, BC1} and
suppose that $f$ and $\sigma $ are independent of $t$.

\begin{proposition}
\label{pr3}Consider 
\begin{equation}
\left\{ 
\begin{array}{lllll}
f^{0} & : & \mathbb{R}^{n} & \mathbb{\rightarrow } & \mathbb{R}^{n}, \\ 
\sigma ^{0} & : & \mathbb{R}^{n} & \mathbb{\rightarrow } & \mathbb{R}%
^{n\times d}.%
\end{array}%
\right.  \tag{4.3}
\end{equation}%
Assume that $f^{0},$ $\sigma ^{0}$ satisfy Lipschitz continuous functions,
with sublinear growth$.$

\noindent Consider 
\begin{equation}
\left\{ 
\begin{array}{lllll}
f^{\varepsilon } & : & \mathbb{R}^{n} & \mathbb{\rightarrow } & \mathbb{R}%
^{n}, \\ 
\sigma ^{\varepsilon } & : & \mathbb{R}^{n} & \mathbb{\rightarrow } & 
\mathbb{R}^{n\times d}.%
\end{array}%
\right.  \tag{4.4}
\end{equation}%
Assume that $f^{\varepsilon }$, $\sigma ^{\varepsilon }$ are uniformly
Lipschitz continuous, have sublinear growth and converge uniformly to $%
f:=f^{0},$ $\sigma :=\sigma ^{0},$ respectively. Then, the family $\left\{
X^{\varepsilon ,t,x}\left( \cdot \right) :0<\varepsilon <1\right\} $ of
random variables of the solutions of the following perturbed stochastic
differential equations%
\begin{equation}
\left\{ 
\begin{array}{rcl}
\mbox{\rm d}X^{\varepsilon ,t,x}\left( r\right) & = & f^{\varepsilon }\left(
X^{\varepsilon ,t,x}\left( r\right) \right) \mbox{\rm d}r+\sqrt{\varepsilon }%
\sigma ^{\varepsilon }\left( X^{\varepsilon ,t,x}\left( r\right) \right) %
\mbox{\rm d}W\left( r\right) , \\ 
X^{\varepsilon ,t,x}\left( t\right) & = & x,%
\end{array}%
\right.  \tag{4.5}
\end{equation}%
obey a large deviations principle in $C\left( \left[ t,T\right] ;\mathbb{R}%
^{n}\right) $, the space of the continuous functions $f:\left[ t,T\right]
\rightarrow \mathbb{R}^{n},$ with the good rate function $\mathcal{I}$
defined as 
\begin{equation*}
\mathcal{I}\left( \phi \right) =\inf \left\{ \frac{1}{2}\left\Vert \varphi
\right\Vert _{\mathcal{H}_{1}\left( t\right) }^{2}:\varphi \in \mathcal{H}%
_{1}\left( t\right) ,\text{ }\phi ^{^{\prime }}\left( s\right) =f\left( \phi
\left( s\right) \right) +\sigma \left( \phi \left( s\right) \right) \varphi
^{\prime }\left( s\right) \right\} .
\end{equation*}%
Moreover, the level sets of $\mathcal{I}$ are compact and for every Borel
subset $A$ of $C\left( \left[ t,T\right] ;\mathbb{R}^{n}\right) ,$ we have 
\begin{equation*}
-\inf\limits_{g\in \mathring{A}}\mathcal{I}\left( g\right) \leq
\liminf\limits_{\varepsilon \rightarrow 0}\varepsilon \log P\left(
X^{\varepsilon ,t,x}\in A\right) \leq \limsup_{\varepsilon \rightarrow
0}\varepsilon \log P\left( X^{\varepsilon ,t,x}\in A\right) \leq
-\inf\limits_{g\in \bar{A}}\mathcal{I}\left( g\right) .
\end{equation*}
\end{proposition}

The proof can be seen in \cite{BC1}. Next we introduce a very important
result in Large deviation theory, used to transfer a LDP from one topology
space to another one.

\begin{lemma}[\textbf{Contraction Principle Theorem 4.2.23 Page 133 in 
\protect\cite{DZ}}]
\label{l9}

Let $\left\{ \mu _{\epsilon }\right\} $ be a family of probability measures
that satisfies the Large Deviation Principle with a good rate function $I$
on a Hausdorff topological space $\mathcal{S}$, and for $\varepsilon \in
\left( 0,1\right] ,$ let $f_{\varepsilon }$: $\mathcal{S}\rightarrow 
\mathcal{Q}$ be continuous functions, with $\left( \mathcal{Q},d\right) $ a
metric space. Assume that there exists a measurable map $f:\mathcal{S}%
\rightarrow \mathcal{Q}$ such that for every $\alpha <+\infty ,$%
\begin{equation*}
\limsup\limits_{\varepsilon \rightarrow 0}\sup\limits_{\left\{ x:I\left(
x\right) \leq \alpha \right\} }d\left( f_{\varepsilon }\left( x\right)
,f\left( x\right) \right) =0.
\end{equation*}%
Then the family of probability measures $\left\{ \mu _{\epsilon }\circ
f_{\varepsilon }^{-1}\right\} $ satisfy the LDP in $\mathcal{Q}$ with the
good rate function 
\begin{equation*}
I^{^{\prime }}\left( y\right) =\inf \left\{ I\left( x\right) :y=f\left(
x\right) \right\} .
\end{equation*}
\end{lemma}

\noindent For our purposes, we give the following

\begin{definition}
Let $t\in \left[ 0,T\right] .$ We define the mapping $F^{\varepsilon
}:C\left( \left[ t,T\right] ;\mathbb{R}^{n}\right) \rightarrow C\left( \left[
t,T\right] ;\mathbb{R}^{n}\right) $ by 
\begin{equation}
F^{\varepsilon }\left( \psi \right) =\left[ s\rightarrow u^{\varepsilon
}\left( s,\psi \left( s\right) \right) \right] ,\quad 0\leq t\leq s\leq T,%
\text{ }\psi \in C\left( \left[ t,T\right] ;\mathbb{R}^{n}\right) , 
\tag{4.6}
\end{equation}%
where $u^{\varepsilon }$ is given by (\ref{2.1}).
\end{definition}

\noindent As we have seen before 
\begin{equation*}
Y^{\varepsilon ,t,x}\left( t\right) =F^{\varepsilon }\left( X^{\varepsilon
,t,x}\right) \left( t\right) .
\end{equation*}%
From now on, for $\varepsilon =0,$ $u$ and $F$ stand for $u^{0}$ and $F^{0}.$

\noindent In order to prove the uniform convergence of the mapping $%
F^{\varepsilon }$, we need the following formula 
\begin{equation*}
\left\Vert F^{\varepsilon }\left( \varphi \right) -F\left( \varphi \right)
\right\Vert =\sup\limits_{ts\in \left[ t,T\right] }\left\vert u^{\varepsilon
}\left( s,\varphi \left( s\right) \right) -u\left( s,\varphi \left( s\right)
\right) \right\vert ,\quad \varphi \in C\left( \left[ t,T\right] ;\mathbb{R}%
^{n}\right)
\end{equation*}

\noindent or 
\begin{equation*}
\left\Vert F^{\varepsilon }\left( \varphi \right) -F\left( \varphi \right)
\right\Vert =\sup\limits_{s\in \left[ t,T\right] }\left\vert Y^{\varepsilon
,s,\varphi \left( s\right) }\left( s\right) -\mathcal{Y}^{s,\varphi \left(
s\right) }\left( s\right) \right\vert .
\end{equation*}

\noindent We have the following:

\begin{theorem}
\label{t3}Assume (A1) and (A2) hold. Then,

i) The family $\left( X^{\varepsilon ,t,x}\left( \cdot \right) \right)
_{\varepsilon \in \left( 0,1\right] }$ satisfy, as $\varepsilon $ goes to $0$%
, a large deviation principle with a rate function%
\begin{equation*}
\mathcal{I}_{1}\left( \phi \right) =\inf_{\left\{ \varphi \in \mathcal{H}%
_{1}\left( t\right) :\text{ }\phi ^{^{\prime }}\left( s\right) =f\left( \phi
\left( s\right) ,u\left( s,\phi \left( s\right) \right) \right) +\sigma
\left( \phi \left( s\right) ,u\left( s,\phi \left( s\right) \right) \right)
\varphi ^{^{\prime }}\left( s\right) ,\text{ }s\in \left[ t,T\right]
\right\} }\frac{1}{2}\left\Vert \varphi \right\Vert _{\mathcal{H}_{1}\left(
t\right) }^{2},
\end{equation*}%
for $g\in C\left( \left[ t,T\right] ;\mathbb{R}^{n}\right) .$

ii) The family $\left( Y^{\varepsilon ,t,x}\left( \cdot \right) \right)
_{\varepsilon \in \left( 0,1\right] }$ satisfy, as $\varepsilon $ goes to $0$%
, a large deviation principle with a rate function%
\begin{equation*}
\mathcal{I}_{2}\left( \phi \right) =\inf \left\{ \mathcal{I}_{1}\left(
\varphi \right) :F\left( \varphi \right) \left( s\right) =\phi \left(
s\right) =u\left( s,\varphi \left( s\right) \right) \text{, }s\in \left[ t,T%
\right] \right\} .
\end{equation*}
\end{theorem}

\proof%
Indeed, noting (A1) and Lemma \ref{l8}, the first assertion has been
obtained in Proposition \ref{pr3}. We are going to prove the second
assertion. By virtue of the contraction principle Lemma \ref{l9}, we just
need to show that $F^{\varepsilon }$, $\varepsilon \in \left(
0,1\right\rceil $ are continuous and $\left\{ F^{\varepsilon }\right\} $
converges uniformly to $F$ on every compact of $C\left( \left[ 0,T\right] ;%
\mathbb{R}^{n}\right) ,$ as $\varepsilon $ tends to zero. \newline
\noindent Proving the \textit{continuity of }$F^{\varepsilon }$\textit{:}

Let $\varepsilon >0$ and $x\in C([t,T],\mathbb{R}^{d})$. Let $(x_{n})_{n\in 
\mathbb{N}}$ be a sequence in $C([t,T],\mathbb{R}^{d})$ converging to $x$ in
the uniform norm. Fix $\delta >0$. Since $\left\Vert x_{n}-x\right\Vert
_{\infty }\rightarrow 0$, there exists $M>0$ such that $\left\Vert
x\right\Vert _{\infty }$ , $\left\Vert x_{n}\right\Vert _{\infty }\leq M$.

Due to Lemma \ref{l1} and Lemma \ref{l7}, we know that $u^\varepsilon $ is a
continuous function in $[0,T] \times \mathbb{R}^d$ and $u^\varepsilon$ is
uniformly continuous in $[t, T] \times K $ where $K = \overline{B(0,M)}
\subset \mathbb{R}^d$.

There exists $\eta >0$ such that for $s,$ $s_{1}\in \lbrack t,T]$ and $z,$ $%
z_{1}\in K,$ $\left\vert s-s_{1}\right\vert <\eta $ and $\left\vert
z-z_{1}\right\vert <\eta $, we have%
\begin{equation*}
\left\vert u^{\varepsilon }(s,z)-u^{\varepsilon }(s_{1},z_{1})\right\vert
<\delta .
\end{equation*}

\noindent Since $x_{n}\rightarrow x$ in $C([t,T],\mathbb{R}^{d})$, fixing $%
n_{0}\in \mathbb{N}$ such that for all $n\geq n_{0}$ we have $\left\Vert
x_{n}-x\right\Vert _{\infty }<\eta $.

For all $r\in \lbrack t,T]$ and for all $n\geq n_{0},$ $x_{n}(r),$ $x(r)\in
K $ and%
\begin{equation*}
\left\vert u^{\varepsilon }(r,x(r))-u^{\varepsilon }(r,x_{n}(r))\right\vert
<\delta .
\end{equation*}

\noindent So we conclude that $F^{\varepsilon }(x_{n})\rightarrow
F^{\varepsilon }(x)$, which proves the continuity of $F^{\varepsilon }$ in
the point $x\in C([t,T],\mathbb{R}^{d})$. Next let us show the uniform
convergence of the mapping $F^{\varepsilon }.$ Consider a compact set $%
\mathcal{K}$ of $C\left( \left[ 0,T\right] ;\mathbb{R}^{n}\right) $ and let 
\begin{equation*}
\mathcal{L}=\left\{ \varphi \left( s\right) ,\text{ }\varphi \in \mathcal{K},%
\text{ }s\in \left[ t,T\right] \right\} .
\end{equation*}%
Clearly, $\mathcal{L}$ is a compact set of $\mathbb{R}^{n}$. By Lemma 6,
there exists a positive constant $\mathcal{C}_{3}$ such that%
\begin{eqnarray*}
&&\sup\limits_{\varphi \in \mathcal{K}}\left\Vert F^{\varepsilon }\left(
\varphi \right) -F\left( \varphi \right) \right\Vert
^{2}=\sup\limits_{\varphi \in \mathcal{K}}\sup\limits_{s\in \left[ t,T\right]
}\left\vert Y^{\varepsilon ,s,\varphi \left( s\right) }\left( s\right) -%
\mathcal{Y}^{s,\varphi \left( s\right) }\left( s\right) \right\vert ^{2} \\
&\leq &\sup\limits_{x\in \mathcal{L}}\sup\limits_{s\in \left[ t,T\right]
}\left\vert Y^{\varepsilon ,s,x}\left( s\right) -\mathcal{Y}^{s,x}\left(
s\right) \right\vert ^{2} \\
&\leq &\mathcal{C}_{3}\sqrt{\varepsilon }.
\end{eqnarray*}%
The proof is complete. 
\endproof%

\begin{remark}
\textsl{Under the conditions (A3) and (A4), we have the same conclusion as
in Theorem \ref{t3}. The proof is similar we omit it.}
\end{remark}

\begin{remark}
\textsl{For fully coupled FBSDEs, that is,}%
\begin{equation*}
\left\{ 
\begin{array}{crl}
X^{\varepsilon ,t,x}\left( s\right) & = & x+\int_{t}^{s}f\left(
X^{\varepsilon ,t,x}\left( r\right) ,Y^{\varepsilon ,t,x}\left( r\right)
,Z^{\varepsilon ,t,x}\left( r\right) \right) \mbox{\rm d}r \\ 
&  & +\sqrt{\varepsilon }\int_{t}^{s}\sigma \left( X^{\varepsilon
,t,x}\left( r\right) ,Y^{\varepsilon ,t,x}\left( r\right) ,Z^{\varepsilon
,t,x}\left( r\right) \right) \mbox{\rm d}W\left( r\right) , \\ 
Y^{\varepsilon ,t,x}\left( s\right) & = & h\left( X^{\varepsilon ,t,x}\left(
T\right) \right) +\int_{s}^{T}g\left( r,X^{\varepsilon ,t,x}\left( r\right)
,Y^{\varepsilon ,t,x}\left( r\right) ,Z^{\varepsilon ,t,x}\left( r\right)
\right) \mbox{\rm d}r \\ 
&  & -\int_{s}^{T}Z^{\varepsilon ,t,x}\left( r\right) \mbox{\rm d}W\left(
r\right) ,\text{ }0\leq t\leq s\leq T,%
\end{array}%
\right.
\end{equation*}%
\textsl{first of all, note that trajectories of the process }$\left(
Z^{\varepsilon ,t,x}\left( s\right) \right) _{\left( t\leq s\leq T\right) }$%
\textsl{\ are not continuous in general. As a matter of fact, under the
assumptions of Theorem 2.6 in \cite{PW}, we know only that }$Z^{\varepsilon
,t,x}$\textsl{\ belongs to the space }$\mathcal{M}^{2}\left( 0,T;\mathbb{R}%
^{m\times d}\right) $\textsl{, which allows us to arbitrarily change the
values of the process }$Z^{\varepsilon ,t,x}$\textsl{\ in any }$P$\textsl{%
-null set. In particular, }$\left\{ t\right\} \times \Omega $\textsl{\ is a }%
$P$\textsl{-null set, which means that }$Z^{\varepsilon ,t,x}$\textsl{\ can
be any }$m\times d$\textsl{\ matrix. Hence, we can not get Lipschitz
property for }$v^{\varepsilon }.$\textsl{\ This issue will be carried out in
our future publications.}
\end{remark}

\appendix

\section{Appendix}

\label{Appendix}

\section*{The Proofs of Lemma \protect\ref{l1}-Lemma \protect\ref{l4}}

\noindent Proof of Lemma \ref{l1}.

\proof
Applying Itô's formula to $\left\vert Y^{\varepsilon ,t,x}\left( \cdot
\right) -Y^{\varepsilon ,t,y}\left( \cdot \right) \right\vert ^{2}$ on $%
\left[ s,T\right] ,$ we have%
\begin{eqnarray}
&&\left\vert Y^{\varepsilon ,t,x}\left( s\right) -Y^{\varepsilon ,t,y}\left(
s\right) \right\vert ^{2}+\mathbb{E}^{\mathcal{F}_{s}}\left[
\int_{s}^{T}\left( \left\vert Z^{\varepsilon ,t,x}\left( r\right)
-Z^{\varepsilon ,t,y}\left( r\right) \right\vert ^{2}\right) \mbox{\rm d}r%
\right]  \notag \\
&=&\mathbb{E}^{\mathcal{F}_{s}}\left[ \left\vert h\left( X^{\varepsilon
,t,x}\left( T\right) \right) -h\left( X^{\varepsilon ,t,y}\left( T\right)
\right) \right\vert ^{2}\right]  \notag \\
&&+2\mathbb{E}^{\mathcal{F}_{s}}\left[ \int_{s}^{T}\left( Y^{\varepsilon
,t,x}\left( r\right) -Y^{\varepsilon ,t,y}\left( r\right) \right) \left(
f^{\varepsilon ,t,x}\left( r\right) -f^{\varepsilon ,t,y}\left( r\right)
\right) \mbox{\rm d}r\right] .  \label{5.1}
\end{eqnarray}%
which yields%
\begin{eqnarray}
&&\mathbb{E}\left[ \sup\limits_{t\leq s\leq T}\left\vert Y^{\varepsilon
,t,x}\left( s\right) -Y^{\varepsilon ,t,y}\left( s\right) \right\vert ^{2}%
\right] +\frac{1}{2}\mathbb{E}\left[ \int_{t}^{T}\left( \left\vert
Z^{\varepsilon ,t,x}\left( r\right) -Z^{\varepsilon ,t,y}\left( r\right)
\right\vert ^{2}\right) \mbox{\rm d}r\right]  \notag \\
&\leq &\mathbb{E}\left[ \left\vert h\left( X^{\varepsilon ,t,x}\left(
T\right) \right) -h\left( X^{\varepsilon ,t,y}\left( T\right) \right)
\right\vert ^{2}\right]  \notag \\
&&+\mathbb{E}\left[ \int_{t}^{T}C_{1}\left\vert X^{\varepsilon ,t,x}\left(
r\right) -X^{\varepsilon ,t,y}\left( r\right) \right\vert ^{2}+\left(
3C_{1}+2C_{1}^{2}\right) \left\vert Y^{\varepsilon ,t,x}\left( r\right)
-Y^{\varepsilon ,t,y}\left( r\right) \right\vert ^{2}\mbox{\rm d}r\right] 
\notag \\
&\leq &C_{1}^{2}\mathbb{E}\left[ \left\vert X^{\varepsilon ,t,x}\left(
T\right) -X^{\varepsilon ,t,y}\left( T\right) \right\vert ^{2}\right]  \notag
\\
&&+\mathbb{E}\left[ \int_{t}^{T}\left( C_{1}\left\vert X^{\varepsilon
,t,x}\left( r\right) -X^{\varepsilon ,t,y}\left( r\right) \right\vert
^{2}+\left( 3C_{1}+2C_{1}^{2}\right) \left\vert Y^{\varepsilon ,t,x}\left(
r\right) -Y^{\varepsilon ,t,y}\left( r\right) \right\vert ^{2}\right) %
\mbox{\rm d}r\right]  \notag \\
&\leq &\left( 3C_{1}+2C_{1}^{2}\right) \Bigg \{\mathbb{E}\left[ \left\vert
X^{\varepsilon ,t,x}\left( T\right) -X^{\varepsilon ,t,y}\left( T\right)
\right\vert ^{2}\right]  \notag \\
&&+\mathbb{E}\left[ \int_{t}^{T}\left( \left\vert X^{\varepsilon ,t,x}\left(
r\right) -X^{\varepsilon ,t,y}\left( r\right) \right\vert ^{2}+\left\vert
Y^{\varepsilon ,t,x}\left( r\right) -Y^{\varepsilon ,t,y}\left( r\right)
\right\vert ^{2}\right) \mbox{\rm d}r\right] \Bigg \},  \label{5.2}
\end{eqnarray}%
since $2ab\leq a^{2}+b^{2}.$

Applying Itô's formula to $X^{\varepsilon ,t,x}\left( s\right)
-X^{\varepsilon ,t,y}\left( s\right) $ yields that%
\begin{eqnarray*}
&&\left\vert X^{\varepsilon ,t,x}\left( s\right) -X^{\varepsilon ,t,y}\left(
s\right) \right\vert ^{2} \\
&=&\int_{t}^{s}2\left( X^{\varepsilon ,t,x}\left( r\right) -X^{\varepsilon
,t,y}\left( r\right) \right) \left( f^{\varepsilon ,t,x}\left( r\right)
-f^{\varepsilon ,t,y}\left( r\right) \right) \mbox{\rm d}r+\varepsilon
\int_{t}^{s}\left\vert \sigma ^{\varepsilon ,t,x}\left( r\right) -\sigma
^{\varepsilon ,t,y}\left( r\right) \right\vert ^{2}\mbox{\rm d}r
\end{eqnarray*}%
\begin{equation}
+\sqrt{\varepsilon }\int_{t}^{s}2\left( X^{\varepsilon ,t,x}\left( r\right)
-X^{\varepsilon ,t,y}\left( r\right) \right) \left( \sigma ^{\varepsilon
,t,x}\left( r\right) -\sigma ^{\varepsilon ,t,y}\left( r\right) \right) %
\mbox{\rm d}W\left( r\right) .  \label{5.3}
\end{equation}%
By Burkholder-Davis-Gundy's inequality, there is a constant $C_{3}>0$ such
that%
\begin{eqnarray}
&&\mathbb{E}\left[ \sup\limits_{t\leq s\leq T}\left\vert X^{\varepsilon
,t,x}\left( s\right) -X^{\varepsilon ,t,y}\left( s\right) \right\vert ^{2}%
\right]  \notag \\
&\leq &C_{3}\Bigg \{\left( 1+\sqrt{\varepsilon }\right) \mathbb{E}\left[
\int_{t}^{T}\left( \left\vert X^{\varepsilon ,t,x}\left( r\right)
-X^{\varepsilon ,t,y}\left( r\right) \right\vert ^{2}\right) \mbox{\rm d}r%
\right] +\mathbb{E}\left[ \int_{t}^{T}\left( \left\vert f^{\varepsilon
,t,x}\left( r\right) -f^{\varepsilon ,t,y}\left( r\right) \right\vert
^{2}\right) \mbox{\rm d}r\right]  \notag \\
&&+\left( \sqrt{\varepsilon }+\varepsilon \right) \mathbb{E}\left[
\int_{t}^{T}\left( \left\vert \sigma ^{\varepsilon ,t,x}\left( r\right)
-\sigma ^{\varepsilon ,t,y}\left( r\right) \right\vert ^{2}\right) 
\mbox{\rm
d}r\right] \Bigg \}  \notag \\
&\leq &C_{3}\mathbb{E}\left[ \int_{t}^{T}\left( \left( 6C_{1}^{2}+2\right)
\left\vert X^{\varepsilon ,t,x}\left( r\right) -X^{\varepsilon ,t,y}\left(
r\right) \right\vert ^{2}+6C_{1}^{2}\left\vert Y^{\varepsilon ,t,x}\left(
r\right) -Y^{\varepsilon ,t,y}\left( r\right) \right\vert ^{2}\right) %
\mbox{\rm d}r\right]  \notag \\
&\leq &C_{3}\mathbb{E}\left[ \left( 6C_{1}^{2}+2\right) \int_{t}^{T}\left(
\left\vert X^{\varepsilon ,t,x}\left( r\right) -X^{\varepsilon ,t,y}\left(
r\right) \right\vert ^{2}\right) \mbox{\rm d}r+6TC_{1}^{2}\sup\limits_{t\leq
s\leq T}\left\vert Y^{\varepsilon ,t,x}\left( s\right) -Y^{\varepsilon
,t,y}\left( s\right) \right\vert ^{2}\right] .  \label{5.4}
\end{eqnarray}%
Applying Itô's formula to $\left( X^{\varepsilon ,t,x}\left( r\right)
-X^{\varepsilon ,t,y}\left( r\right) \right) \left( Y^{\varepsilon
,t,x}\left( r\right) -Y^{\varepsilon ,t,y}\left( r\right) \right) $, and by
virtue of assumption (A2), we have%
\begin{eqnarray}
&&C_{2}\mathbb{E}\left[ \left\vert X^{\varepsilon ,t,x}\left( T\right)
-X^{\varepsilon ,t,y}\left( T\right) \right\vert ^{2}\right] -\mathbb{E}%
\left[ \left( x-y\right) \left( Y^{\varepsilon ,t,x}\left( t\right)
-Y^{\varepsilon ,t,y}\left( t\right) \right) \right]  \notag \\
&\leq &\mathbb{E}\left[ h\left( X^{\varepsilon ,t,x}\left( T\right) \right)
-h\left( X^{\varepsilon ,t,y}\left( T\right) \right) \left( X^{\varepsilon
,t,x}\left( T\right) -X^{\varepsilon ,t,y}\left( T\right) \right) \right] 
\notag \\
&&-\mathbb{E}\left[ \left( x-y\right) \left( Y^{\varepsilon ,t,x}\left(
t\right) -Y^{\varepsilon ,t,y}\left( t\right) \right) \right]  \notag \\
&\leq &-C_{2}\left( 1+\sqrt{\varepsilon }\right) \mathbb{E}\left[
\int_{t}^{T}\left( \left\vert X^{\varepsilon ,t,x}\left( r\right)
-X^{\varepsilon ,t,y}\left( r\right) \right\vert ^{2}+\left\vert
Y^{\varepsilon ,t,x}\left( r\right) -Y^{\varepsilon ,t,y}\left( r\right)
\right\vert ^{2}\right) \mbox{\rm d}r\right]  \notag \\
&\leq &-C_{2}\mathbb{E}\left[ \int_{t}^{T}\left( \left\vert X^{\varepsilon
,t,x}\left( r\right) -X^{\varepsilon ,t,y}\left( r\right) \right\vert
^{2}+\left\vert Y^{\varepsilon ,t,x}\left( r\right) -Y^{\varepsilon
,t,y}\left( r\right) \right\vert ^{2}\right) \mbox{\rm d}r\right] .
\label{5.5}
\end{eqnarray}%
The inequality (\ref{5.5}) can be rewritten as follows%
\begin{eqnarray}
&&C_{2}\mathbb{E}\left[ \left\vert X^{\varepsilon ,t,x}\left( T\right)
-X^{\varepsilon ,t,y}\left( T\right) \right\vert ^{2}\right]  \notag \\
&&+C_{2}\mathbb{E}\left[ \int_{t}^{T}\left\vert X^{\varepsilon ,t,x}\left(
r\right) -X^{\varepsilon ,t,y}\left( r\right) \right\vert ^{2}+\left\vert
Y^{\varepsilon ,t,x}\left( r\right) -Y^{\varepsilon ,t,y}\left( r\right)
\right\vert ^{2}\right] \mbox{\rm d}r  \notag \\
&\leq &\mathbb{E}\left[ \left( x-y\right) \left( Y^{\varepsilon ,t,x}\left(
t\right) -Y^{\varepsilon ,t,y}\left( t\right) \right) \right] .  \label{5.6}
\end{eqnarray}%
Combining (\ref{5.5}) and (\ref{5.6}), we have%
\begin{eqnarray}
&&\mathbb{E}\left[ \sup\limits_{t\leq s\leq T}\left\vert Y^{\varepsilon
,t,x}\left( s\right) -Y^{\varepsilon ,t,y}\left( s\right) \right\vert ^{2}%
\right] +\frac{1}{2}\mathbb{E}\left[ \int_{t}^{T}\left\vert Z^{\varepsilon
,t,x}\left( r\right) -Z^{\varepsilon ,t,y}\left( r\right) \right\vert ^{2}%
\mbox{\rm d}r\right]  \notag \\
&\leq &\frac{\left( 3C_{1}+2C_{1}^{2}\right) }{C_{2}}\mathbb{E}\left[ \left(
x-y\right) \left( Y^{\varepsilon ,t,x}\left( t\right) -Y^{\varepsilon
,t,y}\left( t\right) \right) \right]  \notag \\
&\leq &\frac{1}{2}\left\vert Y^{\varepsilon ,t,x}\left( t\right)
-Y^{\varepsilon ,t,y}\left( t\right) \right\vert ^{2}+\frac{1}{2}\left( 
\frac{3C_{1}+2C_{1}^{2}}{C_{2}}\right) ^{2}\left\vert x-y\right\vert ^{2}.
\label{5.7}
\end{eqnarray}%
Therefore, 
\begin{equation}
\mathbb{E}\left[ \sup\limits_{t\leq s\leq T}\left\vert Y^{\varepsilon
,t,x}\left( s\right) -Y^{\varepsilon ,t,y}\left( s\right) \right\vert ^{2}%
\right] \leq \left( \frac{3C_{1}+2C_{1}^{2}}{C_{2}}\right) ^{2}\left\vert
x-y\right\vert ^{2}.  \label{5.8}
\end{equation}%
By Gronwall's inequality, there exists a positive constant $C_{4},$ such
that 
\begin{equation}
\mathbb{E}\left[ \sup\limits_{t\leq s\leq T}\left\vert X^{\varepsilon
,t,x}\left( s\right) -X^{\varepsilon ,t,y}\left( s\right) \right\vert ^{2}%
\right] \leq C_{4}\left\vert x-y\right\vert ^{2},  \label{5.9}
\end{equation}%
where $C_{4}$ depends on $C_{1},$ $C_{2},$ $C_{3},$ and $T.$ Finally, taking 
$\mathcal{C}_{1}=\max \left\{ \left( \frac{3C_{1}+2C_{1}^{2}}{C_{2}}\right)
^{2},C_{4}\right\} ,$ we get the desired result. 
\endproof%

\noindent Proof of Lemma \ref{l2}

\proof
Applying Itô's formula to $\left\vert X^{\varepsilon ,t,x}\left( s\right)
\right\vert ^{2}$ yields that%
\begin{eqnarray}
\left\vert X^{\varepsilon ,t,x}\left( s\right) \right\vert ^{2}
&=&2\int_{t}^{s}\left( X^{\varepsilon ,t,x}\left( r\right) f^{\varepsilon
,t,x}\left( r\right) \right) \mbox{\rm d}r+\varepsilon
\int_{t}^{s}\left\vert \sigma ^{\varepsilon ,t,x}\left( r\right) \right\vert
^{2}\mbox{\rm d}r  \notag \\
&&+\sqrt{\varepsilon }\int_{t}^{s}2X^{\varepsilon ,t,x}\left( r\right)
\sigma ^{\varepsilon ,t,x}\left( r\right) \mbox{\rm d}B\left( r\right) .
\label{5.10}
\end{eqnarray}%
By Burkholder-Davis-Gundy's inequality, there is a constant $C_{5}>0$ such
that%
\begin{eqnarray}
&&\mathbb{E}\left[ \sup\limits_{t\leq s\leq T}\left\vert X^{\varepsilon
,t,x}\left( s\right) \right\vert ^{2}\right]  \notag \\
&\leq &C_{5}\Bigg \{\mathbb{E}\left[ \int_{t}^{T}\left( \left\vert
X^{\varepsilon ,t,x}\left( r\right) \right\vert ^{2}\right) \mbox{\rm d}%
r+\int_{t}^{T}\left( \left\vert f^{\varepsilon ,t,x}\left( r\right)
\right\vert ^{2}\right) \mbox{\rm d}r+\varepsilon \mathbb{E}\left[
\int_{t}^{T}\left( \left\vert \sigma ^{\varepsilon ,t,x}\left( r\right)
\right\vert ^{2}\right) \mbox{\rm d}r\right] \right] \Bigg \}  \notag \\
&=&C_{5}\Bigg \{\mathbb{E}\bigg [\int_{t}^{T}\left\vert X^{\varepsilon
,t,x}\left( r\right) \right\vert ^{2}+\left\vert f^{\varepsilon ,t,x}\left(
r\right) -f\left( r,0,0\right) +f\left( r,0,0\right) \right\vert ^{2}  \notag
\\
&&+\varepsilon \left\vert \sigma ^{\varepsilon ,t,x}\left( r\right) -\sigma
\left( r,0,0\right) +\sigma \left( r,0,0\right) \right\vert ^{2}\mbox{\rm d}r%
\bigg ]\Bigg \}  \notag \\
&\leq &C_{5}\Bigg \{\mathbb{E}\int_{t}^{T}\left( \left\vert X^{\varepsilon
,t,x}\left( r\right) \right\vert ^{2}+6C_{1}^{2}\left( \left\vert
X^{\varepsilon ,t,x}\left( r\right) \right\vert ^{2}+\left\vert
Y^{\varepsilon ,t,x}\left( r\right) \right\vert ^{2}\right) +2\left\vert
f\left( r,0,0\right) \right\vert ^{2}\right) \text{d}r  \notag \\
&&+4\varepsilon C_{1}^{2}\mathbb{E}\left[ \int_{t}^{T}\left( \left\vert
X^{\varepsilon ,t,x}\left( r\right) \right\vert ^{2}+\left\vert
Y^{\varepsilon ,t,x}\left( r\right) \right\vert ^{2}\right) \text{d}r\right]
+2\int_{t}^{T}\left( \left\vert \sigma \left( r,0,0\right) \right\vert
^{2}\right) \text{d}r\Bigg \}  \notag \\
&\leq &C_{5}\left( 2+8C_{1}^{2}\right) \left[ \mathbb{E}\int_{t}^{T}\left(
\left\vert X^{\varepsilon ,t,x}\left( r\right) \right\vert ^{2}+\left\vert
Y^{\varepsilon ,t,x}\left( r\right) \right\vert ^{2}+\left\vert f\left(
r,0,0\right) \right\vert ^{2}+\left\vert \sigma \left( r,0,0\right)
\right\vert ^{2}\right) \text{d}r\right]  \notag \\
&\leq &C_{5}\left( 2+8C_{1}^{2}\right) \mathbb{E}\bigg [\int_{t}^{T}\left%
\vert X^{\varepsilon ,t,x}\left( r\right) \right\vert ^{2}\text{d}%
r+T\sup\limits_{t\leq s\leq T}\left\vert Y^{\varepsilon ,t,x}\left( s\right)
\right\vert ^{2}  \notag \\
&&+\int_{0}^{T}\left( \left\vert f\left( r,0,0\right) \right\vert
^{2}+\left\vert \sigma \left( r,0,0\right) \right\vert ^{2}\right) \text{d}r%
\bigg ].  \label{5.11}
\end{eqnarray}%
Second applying Itô's formula to $\left\vert Y^{\varepsilon ,t,x}\left(
\cdot \right) \right\vert ^{2}$ on $\left[ s,T\right] ,$ we have%
\begin{eqnarray}
&&\mathbb{E}\left\vert Y^{\varepsilon ,t,x}\left( s\right) \right\vert ^{2}+%
\mathbb{E}^{\mathcal{F}_{s}}\left[ \int_{s}^{T}\left( \left\vert
Z^{\varepsilon ,t,x}\left( r\right) \right\vert ^{2}\right) \text{d}r\right]
\notag \\
&=&\mathbb{E}^{\mathcal{F}_{s}}\left[ \left\vert h\left( X^{\varepsilon
,t,x}\left( T\right) \right) \right\vert ^{2}\right] +2\mathbb{E}^{\mathcal{F%
}_{s}}\left[ \int_{s}^{T}Y^{\varepsilon ,t,x}\left( r\right) f^{\varepsilon
,t,x}\left( r\right) \text{d}r\right] ,  \label{5.12}
\end{eqnarray}%
which yields 
\begin{eqnarray}
&&\mathbb{E}\left[ \sup\limits_{t\leq s\leq T}\left\vert Y^{\varepsilon
,t,x}\left( s\right) \right\vert ^{2}\right] +\frac{1}{2}\mathbb{E}\left[
\int_{t}^{T}\left( \left\vert Z^{\varepsilon ,t,x}\left( r\right)
\right\vert ^{2}\right) \text{d}r\right]  \notag \\
&\leq &\mathbb{E}\left[ \left\vert h\left( X^{\varepsilon ,t,x}\left(
T\right) \right) \right\vert ^{2}+\int_{t}^{T}\left( \left\vert
Y^{\varepsilon ,t,x}\left( r\right) \right\vert ^{2}+\left\vert
f^{\varepsilon ,t,x}\left( r\right) \right\vert ^{2}\right) \text{d}r\right]
\notag \\
&\leq &\mathbb{E}\left[ \left\vert h\left( X^{\varepsilon ,t,x}\left(
T\right) \right) -h\left( 0\right) +h\left( 0\right) \right\vert ^{2}\right]
\notag \\
&&+\mathbb{E}\left[ \int_{t}^{T}\left( \left\vert Y^{\varepsilon ,t,x}\left(
r\right) \right\vert ^{2}+\left\vert f^{\varepsilon ,t,x}\left( r\right)
-f\left( r,0,0,0\right) +f\left( r,0,0,0\right) \right\vert ^{2}\right) 
\text{d}r\right]  \notag \\
&\leq &\mathbb{E}\left[ 2C_{1}^{2}\left\vert X^{\varepsilon ,t,x}\left(
T\right) \right\vert ^{2}+2\left\vert h\left( 0\right) \right\vert ^{2}%
\right]  \notag \\
&&+\mathbb{E}\left[ \int_{t}^{T}\left( \left( \frac{1}{2}+3C_{1}^{2}\right)
\left\vert Y^{\varepsilon ,t,x}\left( r\right) \right\vert ^{2}+\frac{1}{2}%
\left\vert X^{\varepsilon ,t,x}\left( r\right) \right\vert ^{2}+\frac{%
\left\vert f\left( r,0,0,0\right) \right\vert ^{2}}{6C_{1}^{2}}\right) \text{%
d}r\right]  \notag \\
&\leq &\left( 3C_{1}^{2}+3+\frac{1}{6C_{1}^{2}}\right) \bigg \{\mathbb{E}%
\left[ \left\vert X^{\varepsilon ,t,x}\left( T\right) \right\vert
^{2}+\int_{t}^{T}\left( \left\vert Y^{\varepsilon ,t,x}\left( r\right)
\right\vert ^{2}+\left\vert X^{\varepsilon ,t,x}\left( r\right) \right\vert
^{2}\right) \text{d}r\right]  \notag \\
&&+\int_{0}^{T}\left( \left\vert f\left( r,0,0,0\right) \right\vert
^{2}\right) \text{d}r+\left\vert h\left( 0\right) \right\vert ^{2}\bigg \}.
\label{5.13}
\end{eqnarray}%
Set 
\begin{equation*}
C_{6}=\left( 3C_{1}^{2}+3+\frac{1}{6C_{1}^{2}}\right) \left( \left\vert
h\left( 0\right) \right\vert ^{2}+\int_{0}^{T}\left( \left\vert g\left(
r,0,0,0\right) \right\vert ^{2}\right) \text{d}r\right) .
\end{equation*}%
Applying Itô's formula to $X^{\varepsilon ,t,x}\left( r\right)
Y^{\varepsilon ,t,x}\left( r\right) $, and by virtue of assumption (A2), we
have%
\begin{eqnarray}
&&\mathbb{E}\left[ X^{\varepsilon ,t,x}\left( T\right) h^{\varepsilon
,t,x}\left( X^{\varepsilon ,t,x}\left( T\right) \right) \right] -\mathbb{E}%
\left[ xY^{\varepsilon ,t,x}\left( t\right) \right]  \notag \\
&\leq &-C_{2}\left( 1+\sqrt{\varepsilon }\right) \mathbb{E}\left[
\int_{t}^{T}\left( \left\vert X^{\varepsilon ,t,x}\left( r\right)
\right\vert ^{2}+\left\vert Y^{\varepsilon ,t,x}\left( r\right) \right\vert
^{2}\right) \text{d}r\right]  \notag \\
&&+\mathbb{E}\left[ \int_{t}^{T}\left( X^{\varepsilon ,t,x}\left( r\right)
g\left( r,0,0,0\right) +Y^{\varepsilon ,t,x}\left( r\right) f\left(
r,0,0\right) +\sqrt{\varepsilon }Z^{\varepsilon ,t,x}\left( r\right) \sigma
\left( r,0,0\right) \right) \text{d}r\right]  \notag \\
&\leq &-\frac{C_{2}}{2}\mathbb{E}\left[ \int_{t}^{T}\left( \left\vert
X^{\varepsilon ,t,x}\left( r\right) \right\vert ^{2}+\left\vert
Y^{\varepsilon ,t,x}\left( r\right) \right\vert ^{2}\right) \text{d}r\right]
\notag \\
&&+\frac{1}{2C_{2}}\mathbb{E}\left[ \int_{0}^{T}\left( \left\vert g\left(
r,0,0,0\right) \right\vert ^{2}+\left\vert f\left( r,0,0\right) \right\vert
^{2}\right) \text{d}r\right] +\mathbb{E}\left[ \int_{t}^{T}\sqrt{\varepsilon 
}Z^{\varepsilon ,t,x}\left( r\right) \sigma \left( r,0,0\right) \text{d}r%
\right]  \notag \\
&\leq &-\frac{C_{2}}{2}\mathbb{E}\left[ \int_{t}^{T}\left( \left\vert
X^{\varepsilon ,t,x}\left( r\right) \right\vert ^{2}+\left\vert
Y^{\varepsilon ,t,x}\left( r\right) \right\vert ^{2}\right) \text{d}r\right]
+M_{1},  \label{5.14}
\end{eqnarray}%
where%
\begin{equation*}
M_{1}=\frac{\int_{0}^{T}\left\vert g\left( r,0,0,0\right) \right\vert
^{2}+\left\vert f\left( r,0,0\right) \right\vert ^{2}\text{d}r}{2C_{2}}+%
\mathbb{E}\left[ \int_{t}^{T}\left( \sqrt{\varepsilon }Z^{\varepsilon
,t,x}\left( r\right) \sigma \left( r,0,0\right) \right) \text{d}r\right] .
\end{equation*}%
On the other hand, 
\begin{eqnarray}
&&\mathbb{E}\left[ X^{\varepsilon ,t,x}\left( T\right) h^{\varepsilon
,t,x}\left( X^{\varepsilon ,t,x}\left( T\right) \right) \right]  \notag \\
&=&\mathbb{E}\left[ \left( X^{\varepsilon ,t,x}\left( T\right) -0\right)
\left( h^{\varepsilon ,t,x}\left( X^{\varepsilon ,t,x}\left( T\right)
\right) -h^{\varepsilon ,t,x}\left( 0\right) +h^{\varepsilon ,t,x}\left(
0\right) \right) \right]  \notag \\
&\geq &\mathbb{E}\left[ C_{2}\left\vert X^{\varepsilon ,t,x}\left( T\right)
\right\vert ^{2}+X^{\varepsilon ,t,x}\left( T\right) h^{\varepsilon
,t,x}\left( 0\right) \right]  \notag \\
&\geq &\mathbb{E}\left[ C_{2}\left\vert X^{\varepsilon ,t,x}\left( T\right)
\right\vert ^{2}-\frac{\left\vert X^{\varepsilon ,t,x}\left( T\right)
\right\vert ^{2}}{2\alpha }-\alpha \frac{\left\vert h^{\varepsilon
,t,x}\left( 0\right) \right\vert ^{2}}{2}\right] ,  \label{5.15}
\end{eqnarray}%
where $\alpha >0$ large enough such that $C_{2}-\frac{1}{2\alpha }>0.$

Then, we have%
\begin{eqnarray}
&&\left( C_{2}-\frac{1}{2\alpha }\right) \mathbb{E}\left[ \left\vert
X^{\varepsilon ,t,x}\left( T\right) \right\vert ^{2}\right] +\frac{C_{2}}{2}%
\mathbb{E}\left[ \int_{t}^{T}\left( \left\vert X^{\varepsilon ,t,x}\left(
r\right) \right\vert ^{2}+\left\vert Y^{\varepsilon ,t,x}\left( r\right)
\right\vert ^{2}\right) \text{d}r\right]  \notag \\
&\leq &\mathbb{E}\left[ xY^{\varepsilon ,t,x}\left( t\right) \right]
+M_{1}+\alpha \frac{\left\vert h^{\varepsilon ,t,x}\left( 0\right)
\right\vert ^{2}}{2}.  \label{5.16}
\end{eqnarray}%
Setting 
\begin{equation*}
\left\{ 
\begin{array}{l}
M_{2}=M_{1}+\alpha \frac{\left\vert h^{\varepsilon ,t,x}\left( 0\right)
\right\vert ^{2}}{2}, \\ 
\tilde{C}=\min \left\{ C_{2}-\frac{1}{2\alpha },\frac{C_{2}}{2}\right\} ,%
\end{array}%
\right.
\end{equation*}%
we have%
\begin{eqnarray}
&&\mathbb{E}\left[ \left\vert X^{\varepsilon ,t,x}\left( T\right)
\right\vert ^{2}+\int_{t}^{T}\left( \left\vert X^{\varepsilon ,t,x}\left(
r\right) \right\vert ^{2}+\left\vert Y^{\varepsilon ,t,x}\left( r\right)
\right\vert ^{2}\right) \text{d}r\right]  \notag \\
&\leq &\frac{\mathbb{E}\left[ xY^{\varepsilon ,t,x}\left( t\right) \right]
+M_{2}}{\tilde{C}}.  \label{5.17}
\end{eqnarray}%
Noting (\ref{5.13}), we obtain that 
\begin{eqnarray}
&&\mathbb{E}\left[ \sup\limits_{t\leq s\leq T}\left\vert Y^{\varepsilon
,t,x}\left( s\right) \right\vert ^{2}+\frac{1}{2}\int_{t}^{T}\left\vert
Z^{\varepsilon ,t,x}\left( r\right) \right\vert ^{2}\text{d}r\right]  \notag
\\
&\leq &\left( 3C_{1}^{2}+3+\frac{1}{6C_{1}^{2}}\right) \left( \frac{\mathbb{E%
}\left[ xY^{\varepsilon ,t,x}\left( t\right) \right] +M_{2}}{\tilde{C}}%
\right) +C_{6}  \notag \\
&\leq &\frac{\left\vert Y^{\varepsilon ,t,x}\left( t\right) \right\vert ^{2}%
}{2}+\frac{\left\vert x\right\vert ^{2}\left( 3C_{1}^{2}+3+\frac{1}{%
6C_{1}^{2}}\right) ^{2}}{2\tilde{C}^{2}}  \notag \\
&&+\frac{\left( 3C_{1}^{2}+3+\frac{1}{6C_{1}^{2}}\right) M_{2}}{\tilde{C}}%
+C_{6}.  \label{5.18}
\end{eqnarray}%
Define $\tilde{M}=\frac{\left( 3C_{1}^{2}+3+\frac{1}{6C_{1}^{2}}\right) }{%
\tilde{C}}.$ The expression (\ref{5.18}) can be rewritten as%
\begin{eqnarray*}
&&\mathbb{E}\left[ \sup\limits_{t\leq s\leq T}\left\vert Y^{\varepsilon
,t,x}\left( s\right) \right\vert ^{2}+\frac{1}{2}\int_{t}^{T}\left\vert
Z^{\varepsilon ,t,x}\left( r\right) \right\vert ^{2}\text{d}r\right] \\
&\leq &\frac{\left\vert Y^{\varepsilon ,t,x}\left( t\right) \right\vert ^{2}%
}{2}+\frac{\left\vert x\right\vert ^{2}\tilde{M}^{2}}{2}
\end{eqnarray*}
\begin{eqnarray}
&&+\tilde{M}\left( \alpha \frac{\left\vert h^{\varepsilon ,t,x}\left(
0\right) \right\vert ^{2}}{2}+\frac{\int_{0}^{T}\left( \left\vert g\left(
r,0,0,0\right) \right\vert ^{2}+\left\vert f\left( r,0,0\right) \right\vert
^{2}\right) \text{d}r}{2C_{2}}\right)  \notag \\
&&+\tilde{M}^{2}\int_{0}^{T}\left\vert \sigma \left( r,0,0\right)
\right\vert ^{2}\text{d}r+\frac{1}{4}\mathbb{E}\left[ \int_{t}^{T}\left\vert
Z^{\varepsilon ,t,x}\left( r\right) \right\vert ^{2}\text{d}r\right] +C_{6}.
\label{5.19}
\end{eqnarray}%
Consequently, we get%
\begin{eqnarray}
&&\mathbb{E}\left[ \sup\limits_{t\leq s\leq T}\left\vert Y^{\varepsilon
,t,x}\left( s\right) \right\vert ^{2}\right] +\frac{1}{2}\mathbb{E}\left[
\int_{t}^{T}\left\vert Z^{\varepsilon ,t,x}\left( r\right) \right\vert ^{2}%
\text{d}r\right]  \notag \\
&\leq &\left\vert x\right\vert ^{2}\tilde{M}^{2}+\tilde{M}\left( \alpha 
\frac{\left\vert h^{\varepsilon ,t,x}\left( 0\right) \right\vert ^{2}}{2}+%
\frac{\int_{0}^{T}\left( \left\vert g\left( r,0,0,0\right) \right\vert
^{2}+\left\vert f\left( r,0,0\right) \right\vert ^{2}\right) \text{d}r}{%
2C_{2}}\right)  \notag \\
&&+2\tilde{M}^{2}\int_{0}^{T}\left\vert \sigma \left( r,0,0\right)
\right\vert ^{2}\text{d}r+2C_{6}  \notag \\
&\leq &\max \left\{ \Sigma ,\tilde{M}^{2}\right\} \left( 1+\left\vert
x\right\vert ^{2}\right) ,  \label{5.20}
\end{eqnarray}%
where%
\begin{eqnarray*}
\Sigma &=&\tilde{M}\Bigg (\alpha \frac{\left\vert h^{\varepsilon ,t,x}\left(
0\right) \right\vert ^{2}}{2}+\frac{\int_{0}^{T}\left( \left\vert g\left(
r,0,0,0\right) \right\vert ^{2}+\left\vert f\left( r,0,0\right) \right\vert
^{2}\right) \text{d}r}{2C_{2}} \\
&&+2\tilde{M}\int_{0}^{T}\left\vert \sigma \left( r,0,0\right) \right\vert
^{2}\text{d}r\Bigg )+2C_{6}.
\end{eqnarray*}%
By Gronwall's inequality, we derive%
\begin{equation}
\mathbb{E}\left[ \sup\limits_{t\leq s\leq T}\left\vert X^{\varepsilon
,t,x}\left( s\right) \right\vert ^{2}\right] \leq C_{7}\left( 1+\left\vert
x\right\vert ^{2}\right) ,  \label{5.21}
\end{equation}%
where $C_{7}$ is independent of $\varepsilon .$ Taking $\mathcal{C}_{2}=\max
\left\{ \max \left\{ \Sigma ,\tilde{M}^{2}\right\} ,C_{7}\right\} ,$ we get
the desired result. 
\endproof%

\noindent The proof of Lemma \ref{l3}.

\proof
Suppose that $t_{1}>t_{2}.$ By standard estimates and Itô isometry we have%
\begin{eqnarray*}
&&\mathbb{E}\left[ \sup\limits_{t_{1}\leq s\leq T}\left\vert X^{\varepsilon
,t_{1},x}\left( s\right) -X^{\varepsilon ,t_{2},x}\left( s\right)
\right\vert ^{2}\right] \\
&\leq &2\left( t_{1}-t_{2}\right) \mathbb{E}\left[ \int_{t_{2}}^{t_{1}}\left%
\vert f\left( r,X^{\varepsilon ,t_{2},x}\left( r\right) ,Y^{\varepsilon
,t_{2},x}\left( r\right) \right) \right\vert ^{2}\text{d}r\right] \\
&&+2\mathbb{E}\left[ \int_{t_{2}}^{t_{1}}\left\vert \sigma \left(
r,X^{\varepsilon ,t_{2},x}\left( r\right) ,Y^{\varepsilon ,t_{2},x}\left(
r\right) \right) \right\vert ^{2}\text{d}r\right] \\
&\leq &\left( t_{1}-t_{2}\right) \mathbb{E}\Bigg [16TC_{1}^{2}\left(
\sup\limits_{t_{2}\leq r\leq t_{1}}\left\vert X^{\varepsilon ,t_{2},x}\left(
r\right) \right\vert ^{2}+\sup\limits_{t_{2}\leq r\leq t_{1}}\left\vert
Y^{\varepsilon ,t_{2},x}\left( r\right) \right\vert ^{2}\right)
\end{eqnarray*}%
\begin{equation*}
+2\int_{0}^{T}\left\vert f\left( r,0,0\right) \right\vert ^{2}\text{d}%
r+2\int_{0}^{T}\left\vert \sigma \left( r,0,0\right) \right\vert ^{2}\text{d}%
r\Bigg ].
\end{equation*}%
It follows from Lemma \ref{l1} that 
\begin{equation*}
\mathbb{E}\left[ \sup\limits_{t_{1}\leq s\leq T}\left\vert X^{\varepsilon
,t_{1},x}\left( s\right) -X^{\varepsilon ,t_{2},x}\left( s\right)
\right\vert ^{2}\right] \leq \left( t_{1}-t_{2}\right) \Theta ,
\end{equation*}%
where 
\begin{equation*}
\Theta =\left[ 16TC_{1}^{2}\mathcal{C}_{2}\left( 1+\left\vert x\right\vert
^{2}\right) +2\int_{0}^{T}\left\vert f\left( r,0,0\right) \right\vert ^{2}%
\text{d}r+2\int_{0}^{T}\left\vert \sigma \left( r,0,0\right) \right\vert ^{2}%
\text{d}r\right]
\end{equation*}%
Similarly as in Lemma \ref{l1}, we have 
\begin{eqnarray*}
&&\mathbb{E}\left[ \sup\limits_{t_{1}\leq s\leq T}\left\vert Y^{\varepsilon
,t_{1},x}\left( s\right) -Y^{\varepsilon ,t_{2},x}\left( s\right)
\right\vert ^{2}\right] +\frac{1}{2}\mathbb{E}\left[ \int_{t_{1}}^{T}\left%
\vert Z^{\varepsilon ,t_{1},x}\left( r\right) -Z^{\varepsilon
,t_{2},x}\left( r\right) \right\vert ^{2}\text{d}r\right] \\
&\leq &\left( 3C_{1}+2C_{1}^{2}\right) \Bigg \{\mathbb{E}\left[ \left\vert
X^{\varepsilon ,t_{1},x}\left( T\right) -X^{\varepsilon ,t_{2},x}\left(
T\right) \right\vert ^{2}\right] \\
&&+\mathbb{E}\int_{t_{1}}^{T}\left\vert X^{\varepsilon ,t_{1},x}\left(
r\right) -X^{\varepsilon ,t_{2},x}\left( r\right) \right\vert
^{2}+\left\vert Y^{\varepsilon ,t_{1},x}\left( r\right) -Y^{\varepsilon
,t_{2},x}\left( r\right) \right\vert ^{2}\text{d}r\Bigg \}.
\end{eqnarray*}%
By Gronwall's inequality, we obtain 
\begin{equation*}
\mathbb{E}\left[ \sup\limits_{t_{1}\leq s\leq T}\left\vert Y^{\varepsilon
,t_{1},x}\left( s\right) -Y^{\varepsilon ,t_{2},x}\left( s\right)
\right\vert ^{2}\right] \leq \left( t_{1}-t_{2}\right) e^{T}\left(
3C_{1}+2C_{1}^{2}\right) \left( 1+T\right) \Theta ,
\end{equation*}%
and 
\begin{eqnarray*}
&&\mathbb{E}\left[ \int_{t_{1}}^{T}\left\vert Z^{\varepsilon ,t_{1},x}\left(
r\right) -Z^{\varepsilon ,t_{2},x}\left( r\right) \right\vert ^{2}\text{d}r%
\right] \\
&\leq &2\left( t_{1}-t_{2}\right) \left( 3C_{1}+2C_{1}^{2}\right) \Theta
\left( 1+T+e^{T}\left( 3C_{1}+2C_{1}^{2}\right) \left( 1+T\right) \right) .
\end{eqnarray*}%
Set 
\begin{equation*}
\left\{ 
\begin{array}{l}
\mathcal{C}_{5}=e^{T}\left( 3C_{1}+2C_{1}^{2}\right) \left( 1+T\right) \Theta
\\ 
\mathcal{C}_{6}=2\left( 3C_{1}+2C_{1}^{2}\right) \Theta \left(
1+T+e^{T}\left( 3C_{1}+2C_{1}^{2}\right) \left( 1+T\right) \right)%
\end{array}%
\right.
\end{equation*}%
and $\mathcal{C}_{3}=\max \left\{ \Theta ,\mathcal{C}_{5},\mathcal{C}%
_{6}\right\} .$ We get the desired result. 
\endproof%

\noindent The proof of Lemma \ref{l4}.

\proof
Analogously, applying Itô's formula to $\left\vert Y^{\varepsilon
_{1},t,x}\left( \cdot \right) -Y^{\varepsilon _{2},t,x}\left( \cdot \right)
\right\vert ^{2}$ on $\left[ t,T\right] ,$ by the method used above, we have%
\begin{eqnarray*}
&&\mathbb{E}\left[ \sup\limits_{t\leq s\leq T}\left\vert Y^{\varepsilon
_{1},t,x}\left( s\right) -Y^{\varepsilon _{2},t,x}\left( s\right)
\right\vert ^{2}\right] +\frac{1}{2}\mathbb{E}\left[ \int_{t}^{T}\left\vert
Z^{\varepsilon _{1},t,x}\left( r\right) -Z^{\varepsilon _{2},t,x}\left(
r\right) \right\vert ^{2}\text{d}r\right] \\
&\leq &\mathbb{E}\left[ \left\vert h\left( X^{\varepsilon _{1},t,x}\left(
T\right) \right) -h\left( X^{\varepsilon _{2},t,x}\left( T\right) \right)
\right\vert ^{2}\right] \\
&&+\mathbb{E}\left[ \int_{t}^{T}\left( C_{1}\left\vert X^{\varepsilon
_{1},t,x}\left( r\right) -X^{\varepsilon _{2},t,x}\left( r\right)
\right\vert ^{2}+\left( 3C_{1}+2C_{1}^{2}\right) \left\vert Y^{\varepsilon
_{1},t,x}\left( r\right) -Y^{\varepsilon _{2},t,x}\left( r\right)
\right\vert ^{2}\right) \text{d}r\right]
\end{eqnarray*}%
\begin{eqnarray}
&\leq &C_{1}^{2}\mathbb{E}\left[ \left\vert X^{\varepsilon _{1},t,x}\left(
T\right) -X^{\varepsilon _{2},t,x}\left( T\right) \right\vert ^{2}\right] 
\notag \\
&&+\mathbb{E}\left[ \int_{t}^{T}\left( C_{1}\left\vert X^{\varepsilon
_{1},t,x}\left( r\right) -X^{\varepsilon _{2},t,x}\left( r\right)
\right\vert ^{2}+\left( 3C_{1}+2C_{1}^{2}\right) \left\vert Y^{\varepsilon
_{1},t,x}\left( r\right) -Y^{\varepsilon _{2},t,x}\left( r\right)
\right\vert ^{2}\right) \text{d}r\right]  \notag \\
&\leq &\left( 3C_{1}+2C_{1}^{2}\right) \Bigg \{\mathbb{E}\left[ \left\vert
X^{\varepsilon _{1},t,x}\left( T\right) -X^{\varepsilon _{2},t,x}\left(
T\right) \right\vert ^{2}\right]  \notag \\
&&+\mathbb{E}\int_{t}^{T}\left( \left\vert X^{\varepsilon _{1},t,x}\left(
r\right) -X^{\varepsilon _{2},t,x}\left( r\right) \right\vert
^{2}+\left\vert Y^{\varepsilon _{1},t,x}\left( r\right) -Y^{\varepsilon
_{2},t,x}\left( r\right) \right\vert ^{2}\right) \text{d}r\Bigg \}.
\label{5.22}
\end{eqnarray}%
Applying Itô's formula to $X^{\varepsilon _{1},t,x}\left( s\right)
-X^{\varepsilon _{2},t,x}\left( s\right) $ yields that%
\begin{eqnarray}
&&\left\vert X^{\varepsilon _{1},t,x}\left( s\right) -X^{\varepsilon
_{2},t,x}\left( s\right) \right\vert ^{2}  \notag \\
&=&\int_{t}^{s}2\left( X^{\varepsilon _{1},t,x}\left( r\right)
-X^{\varepsilon _{2},t,x}\left( r\right) \right) \left( f^{\varepsilon
_{1},t,x}\left( r\right) -f^{\varepsilon _{2},t,x}\left( r\right) \right) 
\text{d}r+\varepsilon \int_{t}^{s}\left\vert \sigma ^{\varepsilon
_{1},t,x}\left( r\right) -\sigma ^{\varepsilon _{2},t,x}\left( r\right)
\right\vert ^{2}\text{d}r  \notag \\
&&+\sqrt{\varepsilon }\int_{t}^{s}2\left( X^{\varepsilon _{1},t,x}\left(
r\right) -X^{\varepsilon _{2},t,x}\left( r\right) \right) \left( \sigma
^{\varepsilon _{1},t,x}\left( r\right) -\sigma ^{\varepsilon _{2},t,x}\left(
r\right) \right) \text{d}W\left( r\right) .  \label{5.23}
\end{eqnarray}%
By Burkholder-Davis-Gundy's inequality, there is a constant $C_{8}>0$ such
that%
\begin{eqnarray}
&&\mathbb{E}\left[ \sup\limits_{t\leq s\leq T}\left\vert X^{\varepsilon
_{1},t,x}\left( s\right) -X^{\varepsilon _{2},t,x}\left( s\right)
\right\vert ^{2}\right]  \notag \\
&\leq &C_{8}\Bigg \{\mathbb{E}\left[ \int_{t}^{T}\left( 1+\sqrt{\varepsilon }%
\right) \left\vert X^{\varepsilon _{1},t,x}\left( r\right) -X^{\varepsilon
_{2},t,x}\left( r\right) \right\vert ^{2}\text{d}r\right] +\mathbb{E}\left[
\int_{t}^{T}\left\vert f^{\varepsilon _{1},t,x}\left( r\right)
-f^{\varepsilon _{2},t,x}\left( r\right) \right\vert ^{2}\text{d}r\right] 
\notag \\
&&+\left( \sqrt{\varepsilon }+\varepsilon \right) \mathbb{E}\left[
\int_{t}^{T}\left\vert \sigma ^{\varepsilon _{1},t,x}\left( r\right) -\sigma
^{\varepsilon _{2},t,x}\left( r\right) \right\vert ^{2}\text{d}r\right] %
\Bigg \}  \notag \\
&\leq &C_{8}\mathbb{E}\Bigg [\int_{t}^{T}\Big (\left( 6C_{1}^{2}+2\right)
\left\vert X^{\varepsilon _{1},t,x}\left( r\right) -X^{\varepsilon
_{2},t,x}\left( r\right) \right\vert ^{2}+6C_{1}^{2}\left\vert
Y^{\varepsilon _{1},t,x}\left( r\right) -Y^{\varepsilon _{2},t,x}\left(
r\right) \right\vert ^{2}\Big )\text{d}r\Bigg ]  \notag \\
&\leq &C_{8}\mathbb{E}\Bigg [\left( 6C_{1}^{2}+2\right)
\int_{t}^{T}\left\vert X^{\varepsilon _{1},t,x}\left( r\right)
-X^{\varepsilon _{2},t,x}\left( r\right) \right\vert ^{2}\text{d}r  \notag \\
&&+6TC_{1}^{2}\sup\limits_{t\leq s\leq T}\left\vert Y^{\varepsilon
_{1},t,x}\left( s\right) -Y^{\varepsilon _{2},t,x}\left( s\right)
\right\vert ^{2}\Bigg ].  \label{5.24}
\end{eqnarray}%
Once again applying Itô's formula to $\left( X^{\varepsilon _{1},t,x}\left(
r\right) -X^{\varepsilon _{2},t,x}\left( r\right) \right) \left(
Y^{\varepsilon _{1},t,x}\left( r\right) -Y^{\varepsilon _{2},t,x}\left(
r\right) \right) $, and by virtue of assumption (A2), we have%
\begin{eqnarray*}
&&C_{2}\mathbb{E}\left[ \left\vert X^{\varepsilon _{1},t,x}\left( T\right)
-X^{\varepsilon _{2},t,x}\left( T\right) \right\vert ^{2}\right] \\
&\leq &\mathbb{E}\left[ h\left( X^{\varepsilon _{1},t,x}\left( T\right)
\right) -h\left( X^{\varepsilon _{2},t,x}\left( T\right) \right) \left(
X^{\varepsilon _{1},t,x}\left( T\right) -X^{\varepsilon _{2},t,x}\left(
T\right) \right) \right] \\
&\leq &-C_{2}\left( 1+\sqrt{\varepsilon _{1}}\right) \mathbb{E}\left[
\int_{t}^{T}\left( \left\vert X^{\varepsilon _{1},t,x}\left( r\right)
-X^{\varepsilon _{2},t,x}\left( r\right) \right\vert ^{2}+\left\vert
Y^{\varepsilon _{1},t,x}\left( r\right) -Y^{\varepsilon _{2},t,x}\left(
r\right) \right\vert ^{2}\right) \text{d}r\right] \\
&&+\left( \sqrt{\varepsilon _{1}}-\sqrt{\varepsilon _{2}}\right) \mathbb{E}%
\left[ \int_{t}^{T}\sigma ^{\varepsilon _{2},t,x}\left( r\right) \left(
Z^{\varepsilon _{1},t,x}\left( r\right) -Z^{\varepsilon _{2},t,x}\left(
r\right) \right) \text{d}r\right]
\end{eqnarray*}%
\begin{eqnarray}
&\leq &-C_{2}\mathbb{E}\left[ \int_{t}^{T}\left( \left\vert X^{\varepsilon
_{1},t,x}\left( r\right) -X^{\varepsilon _{2},t,x}\left( r\right)
\right\vert ^{2}+\left\vert Y^{\varepsilon _{1},t,x}\left( r\right)
-Y^{\varepsilon _{2},t,x}\left( r\right) \right\vert ^{2}\right) \text{d}r%
\right]  \notag \\
&&+\frac{1}{2}\left( \sqrt{\varepsilon _{1}}-\sqrt{\varepsilon _{2}}\right) 
\mathbb{E}\left[ \int_{t}^{T}\left( \left\vert \sigma ^{\varepsilon
_{2},t,x}\left( r\right) -\sigma \left( r,0,0\right) +\sigma \left(
r,0,0\right) \right\vert ^{2}+\left\vert Z^{\varepsilon _{1},t,x}\left(
r\right) -Z^{\varepsilon _{2},t,x}\left( r\right) \right\vert ^{2}\right) 
\text{d}r\right]  \notag \\
&\leq &-C_{2}\mathbb{E}\left[ \int_{t}^{T}\left( \left\vert X^{\varepsilon
_{1},t,x}\left( r\right) -X^{\varepsilon _{2},t,x}\left( r\right)
\right\vert ^{2}+\left\vert Y^{\varepsilon _{1},t,x}\left( r\right)
-Y^{\varepsilon _{2},t,x}\left( r\right) \right\vert ^{2}\right) \text{d}r%
\right]  \notag \\
&&+\left( \sqrt{\varepsilon _{1}}-\sqrt{\varepsilon _{2}}\right) \Bigg \{%
\mathbb{E}\left[ 4TC_{1}^{2}\left( \sup\limits_{t\leq r\leq T}\left\vert
X^{\varepsilon _{2},t,x}\left( r\right) \right\vert ^{2}+\sup\limits_{t\leq
r\leq T}\left\vert Y^{\varepsilon _{2},t,x}\left( r\right) \right\vert
^{2}\right) \right]  \notag \\
&&+\mathbb{E}\left[ \int_{t}^{T}\left( \left\vert Z^{\varepsilon
_{1},t,x}\left( r\right) \right\vert ^{2}+\left\vert Z^{\varepsilon
_{2},t,x}\left( r\right) \right\vert ^{2}\right) \text{d}r\right]
+2\int_{0}^{T}\left\vert \sigma \left( r,0,0\right) \right\vert ^{2}\text{d}r%
\Bigg \}.  \label{5.25}
\end{eqnarray}%
From Lemma \ref{l2}, we know that 
\begin{eqnarray*}
&&\mathbb{E}\Bigg [\sup\limits_{t\leq r\leq T}\left\vert X^{\varepsilon
_{2},t,x}\left( r\right) \right\vert ^{2}+\sup\limits_{t\leq r\leq
T}\left\vert Y^{\varepsilon _{2},t,x}\left( r\right) \right\vert ^{2} \\
&&+\int_{t}^{T}\left\vert Z^{\varepsilon _{1},t,x}\left( r\right)
\right\vert ^{2}\left( r\right) +\int_{t}^{T}\left\vert Z^{\varepsilon
_{2},t,x}\left( r\right) \right\vert ^{2}\left( r\right) \text{d}r\Bigg ] \\
&\leq &4\mathcal{C}_{2}\left( 1+\left\vert x\right\vert ^{2}\right) .
\end{eqnarray*}%
Hence, combining (\ref{5.22}), (\ref{5.25}), (A1) and (A2), we have%
\begin{eqnarray}
&&\mathbb{E}\left[ \sup\limits_{t\leq s\leq T}\left\vert Y^{\varepsilon
_{1},t,x}\left( s\right) -Y^{\varepsilon _{2},t,x}\left( s\right)
\right\vert ^{2}\right] +\frac{1}{2}\mathbb{E}\left[ \int_{t}^{T}\left\vert
Z^{\varepsilon _{1},t,x}\left( r\right) -Z^{\varepsilon _{2},t,x}\left(
r\right) \right\vert ^{2}\text{d}r\right]  \notag \\
&\leq &\frac{\left( \sqrt{\varepsilon _{1}}-\sqrt{\varepsilon _{2}}\right) }{%
2}M_{3}\left( 3C_{1}+2C_{1}^{2}\right) ,  \label{5.26}
\end{eqnarray}%
where 
\begin{equation*}
M_{3}=\frac{\left( 3C_{1}+2C_{1}^{2}\right) \left[ \left( 1+\left\vert
x\right\vert ^{2}\right) \left( 2\mathcal{C}_{2}+8TC_{1}^{2}\mathcal{C}%
_{2}\right) +2\int_{0}^{T}\left\vert \sigma \left( r,0,0\right) \right\vert
^{2}\text{d}r\right] }{C_{2}}.
\end{equation*}%
By Gronwall's inequality, we also obtain 
\begin{equation*}
\mathbb{E}\left[ \sup\limits_{t\leq s\leq T}\left\vert X^{\varepsilon
_{1},t,x}\left( s\right) -X^{\varepsilon _{2},t,x}\left( s\right)
\right\vert ^{2}\right] \leq C_{9}\left( \sqrt{\varepsilon _{1}}-\sqrt{%
\varepsilon _{2}}\right) ,
\end{equation*}%
where $C_{9}$ is independent of $\varepsilon .$ Taking $\mathcal{C}_{4}=\max
\left\{ C_{9},\frac{M_{3}\left( 3C_{1}+2C_{1}^{2}\right) }{2}\right\} ,$ we
complete the proof. 
\endproof%

\textbf{Acknowledgements}. The first and second authors wish to thank
Fernanda Cipriano for reading this work, and for her constructive
suggestions and questions. The third author wishes to thank Prof. Zhen Wu
and Dr Zhiyong Yu for helpful suggestion and conversations and Prof. Fuqing
Gao for providing the reference \cite{Ba}.

The authors are also very grateful to the referee and the AE for their
fruitful comments and suggestions.

\end{document}